\documentclass[10pt]{article}
\usepackage{latexsym}
\usepackage{amsfonts}
\usepackage{amssymb}
\usepackage{amsmath} 
\usepackage[active]{srcltx}
\usepackage{amsthm}
\usepackage{comment}

\newcommand{\be}[1]{\begin{equation} \label{#1}}
\newcommand{\ee}{\end{equation}}
\def\rf#1{(\ref{#1})}

\numberwithin{equation}{section}

\newtheorem{thm}{Theorem}[section]
\newtheorem{lem}[thm]{Lemma}
\newtheorem{cor}[thm]{Corollary}
\newtheorem{df}[thm]{Definition}
\newtheorem{rem}[thm]{Remark}
\newcommand{\bt}[1]{\begin{thm}}
\newcommand{\et}{\end{thm}}
\newcommand{\bl}[1]{\begin{lem} \label{#1}}
\newcommand{\el}{\end{lem}}
\newcommand{\bc}[1]{\begin{cor} \label{#1}}
\newcommand{\ec}{\end{cor}}
\newcommand{\br}[1]{\begin{rem}\rm \label{#1}}
\newcommand{\er}{\end{rem}}

\newcommand{\bd}[1]{\begin{df} \label{#1}}
\newcommand{\ed}{\end{df}}

\newcommand{\wlims}{\stackrel{*}{\rightharpoonup}}

\newcommand{\rga}{\rightarrow}

\newcommand{\intaverage}{\mbox{$
  {\mbox{$\int$}} \hspace{-.22cm} \mbox{\raisebox{0.17ex}{\mbox{--}}}$} }
\newcommand{\linea}{\mathop{\hbox{
\vrule height .5pt width 6pt depth 0pt}}\nolimits}
\newcommand{\intav}{\mbox{$
  {\mbox{$\displaystyle\int$}} \hspace{-.43cm} 
\mbox{\raisebox{0.59ex}{\mbox{$\linea$}}}$} }
\newcommand{\blist}{\begin{list}{(\roman{icount})}{\usecounter{icount}}}
\newcommand{\elist}{\end{list}\setcounter{icount}{1}}

\newcommand{\Om}{\Omega}
\newcommand{\eps}{\varepsilon}

\newcommand{\cM}{{\cal M}}

\newcommand{\cF}{{\cal F}}
\newcommand{\cG}{{\cal G}}
\newcommand{\cQ}{{\cal Q}}

\newcommand{\Ln}{{\cal L}^{N}}

\newcommand{\mres}{\mathbin{\vrule height 1.6ex depth 0pt width
0.13ex\vrule height 0.13ex depth 0pt width 1.3ex}}

\newcommand\restr[2]{{
  \left.\kern-\nulldelimiterspace 
  #1 
  \vphantom{\big|} 
  \right|_{#2} 
  }}

\newcommand{\Rn}{\mathbb{R}^N}       
       
\newcommand{\Rd}{\mathbb{R}^d}      
\newcommand{\Rdn}{\mbox{${\mathbb R}^{d\times N}$}}

\newcommand{\BHO}{BH(\Omega;\Rd)}

\newcommand{\Qeps}{Q(x_0,\eps)}

\newcommand{\sing}{\frac{d D_s(\nabla u)}{d |D_s(\nabla u)|}}

\title{}
\author{{\Large } \\ Department of Mathematics \\
\and {\Large }}

\begin{document}
\newcounter{icount}

\begin{center}
{\Large \bf Relaxation of Functionals in the Space of Vector-Valued Functions of Bounded Hessian}
\vspace{.5cm} \\
{ \bf  Adrian Hagerty } \\
{ \bf Carnegie Mellon University } \\
{ \bf ahagerty@andrew.cmu.edu }
\vspace*{.5cm} \\
\today
\end{center}
\begin{abstract}
In this paper it is shown that if $\Om \subset \mathbb{R}^N$ is an open, bounded Lipschitz set, and if $f: \Om \times \mathbb{R}^{d \times N \times N} \rga [0, \infty)$ is a continuous function with $f(x, \cdot)$ of linear growth for all $x \in \Om$, then the relaxed functional in the space  of functions of Bounded Hessian of the energy
\[ F[u] = \int_{\Om} f(x, \nabla^2u(x)) dx \]
for bounded sequences in $W^{2,1}$ is given by
\[ \cF[u] = \int_\Om \cQ_2f(x, \nabla^2u) dx + \int_\Om (\cQ_2f)^{\infty}\bigg(x, \frac{d D_s(\nabla u)}{d |D_s(\nabla u)|} \bigg) d |D_s(\nabla u) |. \]
This result is obtained using blow-up techniques and establishes a second order version of the $BV$ relaxation theorems of Ambrosio and Dal Maso \cite{ADM92} and Fonseca and M\"uller \cite{FM93}. The use of the blow-up method is intended to facilitate future study of integrands which include lower order terms and applications in the field of second order structured deformations.
\vspace*{.25cm}

\noindent \textit{Mathematics Subject Classification (2010):} 49J45, 49Q20 \\
\noindent \textit{Keywords: Bounded Hessian functionals, quasiconvexity, relaxation.}
\end{abstract}

\noindent\makebox[\linewidth]{\rule{\paperwidth}{0.4pt}}

\section{Introduction}

Variational models involving higher order derivatives of the underlying fields are ubiquitous in the mathematical analysis of materials and imaging science. As remarked by Demengel in \cite{D84}, the space of functions of bounded Hessian, $BH$, is the natural setting for the study of certain problems in plasticity and elasto-plasticity (see e.g., \cite{BCDMP16} \cite{CLT92} \cite{CLT04} \cite{DM13} \cite{DM15} \cite{D83} \cite{MM16} \cite{PT05}). Similar function spaces have been adopted in the study of structured deformations, which model geometrical changes at microscopic and macroscopic scales. Here the first order theory fails to account for the effect of microscopic jumps in the gradients and curvature effects, and the second order theory was introduced in \cite{OP00} using the space $BH$ and related spaces $SBH, SBV^2$. Recent results in \cite{BMMO17} discuss the second order theory in the $SBV^2$ setting and establish relaxation and integral representation theorems.

Recently, the field of image processing has seen progress using energies with second order terms. The addition of a second order term in the Rudin-Osher-Fatemi TV denoising model can act as a regularizing factor, preventing the so-called ``staircasing effect'', as discussed in \cite{BP10}, \cite{BP13},  \cite{B16}, \cite{DFLM09}. More applications of second order terms in regularization and denoising may be found in \cite{FM17}, \cite{HS06}, \cite{S15}. Further, second order energies have seen applications in variational image fusion \cite{LZ16} and image colorization \cite{JK16}.

In this paper we prove an integral representation result in $BH$ for relaxed functionals with linear growth defined on $W^{2,1}$. To be precise, given an open bounded Lipschitz set $\Om$ in $\Rn$ we define the functional

$$
F[u] := \int_\Om f(x, \nabla^2u) dx, \, \, \, \, u \in W^{2,1}(\Om,\Rd)
$$
\noindent
where $f: \Om \times \mathbb{R}^{d \times N \times N} \rga [0, \infty)$ is a continuous function with $f(x, \cdot)$ of linear growth for all $x \in \Om$, we consider the weakly-$*$ lower semicontinuous envelope in the space $BH(\Om, \Rd)$, namely

$$
\cF[u] := \inf \bigg\{ \liminf_{n \rga \infty} F[u_n] : u_n \rga u \textrm{ in } L^1(\Om; \mathbb{R}^d), \ \sup_{n} \|u_n\|_{W^{2,1}} < \infty \bigg\}. 
$$
\noindent
In Theorem \ref{mainresult} we prove that that, under some growth and regularity hypotheses on $f$, we can represent $\cF$ as 

$$
\cF[u] = \int_\Om \cQ_2f(x, \nabla^2u) dx + \int_\Om (\cQ_2f)^{\infty}\bigg(x, \frac{d D_s(\nabla u)}{d |D_s(\nabla u)|} \bigg) d |D_s(\nabla u) |,
$$

\noindent
where $\cQ_2 f$ is the 2-quasiconvex envelope of $f$, and for $x \in \Om$, $(\cQ_2 f)^{\infty}(x, \cdot)$ is the \textit{recession function} of $\cQ_2f(x, \cdot)$, defined via

\be{recdef}
g^{\infty}(H) := \limsup_{t \rga + \infty} \frac{g(tH)}{t},
\ee
for a function $g: \mathbb{R}^{d \times N \times N} \rga [0, \infty)$.

This result may be seen as a second order version of those of Ambrosio and Dal Maso \cite{ADM92} as well as Fonseca and M\"uller \cite{FM93} on first order linear growth functionals.

This paper is inspired by recent progress in the field of $\mathcal{A}$-quasiconvexity, in the sense of Fonseca and M\"uller as introduced in \cite{FM99} (see also \cite{D82}). In \cite{ADPR17}, Arroyo-Rabasa, De Philippis, and Rindler use a combination of the blow-up method and Young measures to prove a relaxation result in a very general setting. Knowing that the space $BH$ may be viewed through the lens of $\mathcal{A}$-quasiconvexity, this paper adopts some of these techniques to $BH$ relaxation, using the explicit structure of $BH$ that leads to deterministic arguments that avoid the use of Young measures.

In a recent relaxation paper, Breit, Diening and Gmeineder \cite{BDG17} examine what they call $\mathbb{A}$-quasiconvexity. As they note in Section 5, with the existence of an annihilator $\mathbb{L}$, this reduces to what Dacorogna called A-B quasiconvexity \cite{D82}, with $A := \mathbb{L}$ and $B := \mathbb{A}$. Although $BH$ may be viewed within the context of $\mathcal{A}$-quasiconvexity by restricting matrix-valued measures to lie in the subspace $\mathbb{R}^d \otimes \mathbb{R}^{N \times N}_{sym}$, it is not obvious that we can cast it in the framework of A-B quasiconvexity as easily. Regardless, for their argument of lower-semicontinuity, Breit, Diening and Gmeineder make use of $\cite{ADPR17}$, and thus do not provide a Young measure-free argument.

The goal of this paper is to establish a relaxation result in the space $BH$ using purely blow-up methods without reference to Young measures. In order to establish the upper bound, we prove an area-strict density theorem in a very general setting, requiring no extra structure on the limiting measure $\mu$. We develop the direct argument with an eye towards ultimately including lower order terms in the relaxation, which is necessary for related problems, including a second order theory of structured deformations which allows for a Cantor part in the Hessian (see \cite{FHP18}).

While the arguments in this paper are only for the second order case, an extension to higher order derivatives should be possible without much extra work. We note that a higher order relaxation result is addressed by Amar and De Cicco in \cite{AdC94}. However, to our knowledge there seems to be a gap in the proof of lower semicontinuty in \cite{AdC94}, in particular in what concerns the singular part. The proof in our paper makes use of modern results which take a completely different approach in proving lower semicontinuity.

This paper is structured as follows: In Section 2 we establish some preliminary Reshetynak continuity results and geometric properties of 2-quasiconvexity. In Section 3 we prove an area-strict density result for Radon measures (see Theorem \ref{density}), and we obtain a second order extension theorem in order to apply this theorem to $BH$. Section 4 contains the main result, Theorem \ref{mainresult}, which is achieved by a direct blow-up argument.

\

\section{Preliminaries}
Let $d, N \in \mathbb{N}$, let $\Om$ be an open, bounded subset of $\Rn$, and consider the space of functions of Bounded Hessian

\begin{eqnarray*}
\mbox{BH}(\Omega;\Rd) & := & \{u\in \mbox{W}^{1,1}(\Omega;\Rd):D^2 u
\mbox{ is a finite Radon measure}\} \\
& = & \{u\in \mbox{L}^1(\Omega;\Rd):Du\in\mbox{BV}(\Omega;\Rdn)\}.
\end{eqnarray*}
\noindent

We define the unit cube $Q := \{ x \in \Rn : |x_i| \leq \frac{1}{2} \textrm{ for all } 1 \leq i \leq N \}$, and we consider the cube of side length $r$ centered at $x_0$,  $Q(x_0, r) := x_0 + rQ = \{x_0 + ry : y \in Q \}. $

\

For simplicity of notation, we will often denote the Lebesgue measure of a Borel set $E \subset \Rn$ via the notation $|E| := \Ln(E)$.

\

We recall that if $\mu \in \cM( \Om ; \Rd)$ is a finite Radon measure, then there exist $\mu_{ac} \in L^1(\Om; \Rd)$, $\mu_s \in \cM(\Om; \Rd)$ with $|\mu_s| \perp \Ln \mres \Om $, and a $|\mu_s|$-measurable function $\nu_{\mu} = \frac{d\mu_s}{d|\mu_s|}$ with $|\nu_{\mu}(x)| = 1$ for $|\mu_s|$ almost every $x \in \Om$, such that for every Borel set $E \subset \Om$ we have

$$\mu(E) = \int_{E} \mu_{ac}(x) dx + \int_{E} \nu_{\mu}(x) d|\mu_s|(x),$$
and we write $\mu = \mu_{ac} \Ln \mres \Om + \mu_s$. In what follows, if the measure $\mu$ being referenced is clear in the context, we will often drop the subscript and write $\nu_{\mu}$ as $\nu$. Similarly, the total variation $|\mu|$ admits a Radon-Nikodym decomposition $|\mu| = |\mu|_{ac} \Ln \mres \Om + |\mu|_s$. The next lemma establishes that the total variation of the absolutely continuous and the singular part of a Radon measure coincide, respectively, with the absolutely continuous and singular part of the total variation. 

\bl{Radon Nikodym}

Let $\mu \in \cM( \Om; \Rd)$ be a finite Radon measure. Then $|\mu|_{ac}(x) =  |\mu_{ac}|(x)$ for $\Ln$ almost every $x \in \Om$ and $|\mu_s| = |\mu|_s$.

\el
\proof

By the definition of a singular measure and the Radon-Nikodym decomposition, we can split $\Om$ into disjoint Borel sets $A,B$ such that $\mu \mres A = \mu_{ac} \Ln \mres \Om$, $\mu \mres B = \mu_s$, $|\mu_s|(A) = |B| = 0$.  We claim that
$$|\mu| = |\mu| \mres A + |\mu| \mres B$$
is the Radon-Nikodym decomposition of $|\mu|$ with respect to $\Ln \mres \Om$. Indeed, for any Borel set $E$ such that $|E| = 0$, we have 
$$(|\mu| \mres A)(E) = |\mu|(E \cap A) = \sup \bigg\{ \sum_{i = 1}^{\infty} |\mu(F_i)| : F_i \textrm{ Borel, pairwise disjoint, } \bigcup_{i = 1}^{\infty} F_i = E \cap A \bigg \}. $$
Note that for any $F \subset E \cap A$, we have $|\mu(F)| = |\mu(F \cap A)| = |(\mu \mres A)(F)| = 0$, since $\mu \mres A << \Ln$ and $F \subset E$. We conclude that $(|\mu| \mres A)(E) = 0$ and thus $|\mu| \mres A << \Ln$. On the other hand, we claim that the measure $|\mu| \mres B$ is singular with respect to the Lebesgue measure. This follows from the fact that $|B| = 0$ and $(|\mu| \mres B)(A) = |\mu|(A \cap B) = 0$ since $A$ and $B$ are disjoint. We conclude $|\mu| \mres A = |\mu|_{ac} \Ln \mres \Om$ and $|\mu| \mres B = |\mu|_s$.

Now, we show that for any Borel set $P$ we have $|\mu| \mres P = |\mu \mres P|$. Indeed, if $E$ is a Borel set, then
\begin{align}
\notag
|\mu \mres P|(E) &= \sup \bigg\{ \sum_{i = 1}^{\infty} |\mu \mres P(F_i)| : F_i \textrm{ Borel, pairwise disjoint, } \bigcup_{i = 1}^{\infty} F_i = E \bigg \} \\
\notag
&= \sup \bigg\{ \sum_{i = 1}^{\infty} |\mu(P \cap F_i)| : F_i \textrm{ Borel, pairwise disjoint, } \bigcup_{i = 1}^{\infty} F_i = E \bigg \}  \\
\label{stars}
&\leq \sup \bigg\{ \sum_{i = 1}^{\infty} |\mu(G_i)| : G_i \textrm{ Borel, pairwise disjoint, } \bigcup_{i = 1}^{\infty} G_i = E \cap P \bigg \} \\
\label{one}
&= |\mu|(E \cap P) = (|\mu| \mres P)(E).
\end{align}
Given any family $\{G_i\}_{i = 1}^{\infty}$ as in (\ref{stars}), we can add $G_0 := E \setminus P$ which is disjoint from each of the $G_i$, satisfies $|\mu(P \cap G_0)| = 0$ and $\bigcup_{i = 0}^{\infty} G_i = E$. Thus, we can bound (\ref{stars})  from above by
\[ \sup \bigg\{ \sum_{i = 0}^{\infty} |\mu(P \cap G_i)| : G_i \textrm{ Borel, pairwise disjoint, } \bigcup_{i = 0}^{\infty} G_i = E \bigg \} = |\mu \mres P|(E). \]
This, together with (\ref{one}), yields $|\mu \mres P|(E) = (|\mu| \mres P)(E)$ for every Borel set $E$, and so the measures are equal. Thus we have
$$|\mu_{ac}| \Ln \mres \Om = |\mu_{ac} \Ln \mres \Om| = |\mu \mres A| = |\mu| \mres A = |\mu|_{ac} \Ln \mres \Om,$$
therefore $|\mu_{ac}(x)| = |\mu|_{ac}(x)$ for $\Ln$ almost every $x \in \Om$ and
$$|\mu_s| = |\mu \mres B| = |\mu| \mres B = |\mu|_s.$$ 
\qed

\

We consider a functional $F: W^{2,1}(\Om,\Rn) \rightarrow [0, \infty]$ given by

$$
F[u] := \int_\Om f(x, \nabla^2u) dx, \, \, \, \, u \in W^{2,1}(\Om,\Rn),
$$
\noindent
where $f : \Om \, \times \, \mathbb{R}^{d \times N \times N} \rga [0, \infty)$ is a continuous integrand satisfying the following hypotheses:

\blist
\item[(H1)] Linear growth at infinity: $0 \leq f(x,H) \leq C (1 + |H|)$ for all $x \in \Om$, $H \in \mathbb{R}^{d \times N \times N}$ and some $C>0$; 

\item[(H2)] Modulus of continuity: $|f(x, H) - f(y,H)| \leq \omega(|x-y|)(1 + |H|)$ for all $x,y \in \Om, H \in \mathbb{R}^{d \times N \times N}$, where $\omega(s)$ is a nondecreasing function with $w(s) \rga 0$ as $s \rga 0^+$.

In the proof of the main result, we will begin by demonstrating the lower bound under the additional assumption:

\item[(H3)] Coercivity: $f(x,H) \geq c |H| \textrm{ for all } x \in \Om, H \in \mathbb{R}^{d \times N \times N} \textrm{ and some } c \in (0,1).$


\elist

We say that a Borel measurable function $g : \mathbb{R}^{d \times N \times N} \rga [0, \infty)$ is \textit{2-quasiconvex} if for every $H \in \mathbb{R}^{d \times N \times N}$ and for all $w \in W^{1,2}_0(Q; \Rd)$, we have
\be{defQC}
g(H) \leq \int_{Q} g(H+ \nabla^2 w(z)) dz.
\ee
The notion of 2-quasiconvexity, introduced by Meyers in \cite{M65}, is an extension of quasiconvexity (see \cite{D08}) to second order integrands. In general, $k$-quasiconvexity is known to be a necessary and sufficient condition for lower semicontinuity of integrands of $k$th order derivatives in the Sobolev setting. 

We also recall that the \textit{2-quasiconvex envelope} of $g$ is given by
$$\cQ_2g(H) := \inf \bigg\{ \int_{Q} g(H + \nabla^2 w(z) ) dz \  | \ w \in W^{1,2}_0(Q; \Rd) \bigg \} \textrm{ for all } H \in \mathbb{R}^{d \times N \times N}.$$ 

From \cite{FM99}, Proposition 3.4, we have that for an upper semicontinuous function $g : \mathbb{R}^{d \times N \times N} \rga \mathbb{R}$, $\cQ_2 g$ is 2-quasiconvex. Given a continuous function $f: \Om \times \mathbb{R}^{d \times N \times N} \rga [0, \infty)$, with an abuse of language we say that $f$ is quasiconvex if $f(x, \cdot)$ is quasiconvex for all $x \in \Om$, and $\cQ_2f(x, \cdot)$ denotes the 2-quasiconvex envelope of $f(x, \cdot)$.

\bl{inherit}
Let $f : \Om \times \mathbb{R}^{d \times N \times N} \rga [0, \infty)$ be a continuous function satisfying (H1), (H2) and (H3). Then
$$|\cQ_2f(x, H) - \cQ_2f(y,H)| \leq \frac{2C}{c} \omega(|x-y|)(1 + |H|) \textrm{ for all } x,y \in \Om, H \in \mathbb{R}^{d \times N \times N}$$
where $\omega, c, C$ are as in (H1)-(H3).
\el

\proof
For any $x, y \in \Om$, $H \in \mathbb{R}^{d \times N \times N}$ and $w \in W^{2,1}_0(Q; \Rd)$ we have

\begin{align}
\notag
\int_{Q} f(x, &H + \nabla^2 w(z)) dz - \int_{Q} f(y, H + \nabla^2 w(z)) dz \\
\notag
\leq &\int_{Q} |f(x, H+ \nabla^2 w(z)) - f(y, H+ \nabla^2w(z))| dz\\
\leq &\int_{Q} \omega (|x-y|) (1 + |H+\nabla^2w(z)|) dz 
\label{modcont}
\end{align}
by (H2). If $f$ is coercive, then
\begin{equation}
\int_{Q} c \omega (|x-y|)  (1 + |H+ \nabla^2 w(z)|) dz \leq   \int_{Q} \omega (|x-y|)  (1 + f(y,H+\nabla^2 w(z)) dz,
\label{modcont2}
\end{equation}
and so we have by (\ref{modcont}) and (\ref{modcont2}),
\be{modcont3}
\cQ_2 f(x,H) - \int_{Q} f(y, H + \nabla^2 w(z)) dz \leq  \frac{1}{c} \int_{Q} \omega(|x-y|)  (1 + f(y,H+ \nabla^2 w(z)) dz.
\ee
For $\eta >0$ fixed, choose $w$ such that
$$\int_Q f(y, H + \nabla^2 w(z)) dz \leq \cQ_2f(y,H) + \eta,$$
so that (\ref{modcont3}) yields
$$\cQ_2 f(x,H) - \cQ_2 f(y,H) - \eta \leq \frac{1}{c} \omega(|x-y|) (1 + \cQ_2 f(y,H) + \eta ) \leq 2\frac{C}{c} \omega(|x-y|) (1 + |H| + \eta) $$
where we used the fact that, by (H1), we have $\cQ_2f(y,H) \leq C(1+|H|)$. Letting $\eta \rga 0$ we have
$$\cQ_2 f(x,H) - \cQ_2 f(y,H) \leq \tilde{C} \omega(|x-y|) (1 + |H|) $$
and by symmetry, the inequality holds where $x$ and $y$ are switched, yielding
\be{modcontoff}
|\cQ_2 f(x,H) - \cQ_2 f(y,H)| \leq \tilde{C} \omega(|x-y|) (1 + |H|).
\ee
\qed

We will make use of the following lemma concerning the diagonalization of doubly indexed sequences.

\begin{lem}
\label{diags}
Let $X$ be a separable metric space, and let $\{\mu_{n,k}\}, \{\mu_n\} \subset \cM(X; \Rd)$, $\mu \in \cM(X;V)$ be such that
$$ \mu_{n,k} \wlims \mu_n \textrm{ as } k \rga \infty \textrm{ for every } n \in \mathbb{N}, $$
$$ \mu_n \wlims \mu, $$
and
$$\limsup_{n \rga \infty} \limsup_{k \rga \infty} \|\mu_{n,k}\| < \infty. $$
Then, for every $n \in \mathbb{N}$ there exists $k_n$ such that
$$\mu_{n,k_n} \wlims \mu \textrm{ as } n \rga \infty. $$
\end{lem}

The proof of this lemma uses standard density arguments and we omit it here.

To prove Theorem \ref{density} we will apply a lower semicontinuity theorem of Reshetnyak (see \cite{S11}).

\bt{reshet}
\label{reshet}
Let $X$ be a locally compact, separable metric space and let $\{ \mu_n \}$ be a sequence in $\cM(X; \Rd)$. If $\mu_n \wlims \mu \in \cM(X; \Rd)$, then

$$\liminf_{n \rga \infty} \int_X G\bigg( \frac{d \mu_n}{d|\mu_n|}(x) \bigg) d|\mu_n| \geq \int_X G\bigg(\frac{d \mu}{d|\mu|}(x) \bigg) d|\mu| $$ 
\noindent
for every continuous $G : \Rd \rga \mathbb{R}$ positively 1-homogeneous and convex, satisfying the growth condition $|G(\xi)| \leq C|\xi|$ for each $\xi \in \Rd$ and for some $C >0$.
\et

Our relaxation result relies on geometric properties of Hessians and 2-quasiconvex functions. In particular, any 2-quasiconvex function is convex along certain directions, in analogy with quasiconvex functions being rank-one convex (\cite{D08}, Theorem 1.7). More precisely, we will follow the notation of Ball, Currie, and Olver \cite{BCO81}. By $X(N,d,2)$ we denote the space of symmetric bilinear maps from $\mathbb{R}^N \times \mathbb{R}^N$ into $\mathbb{R}^d$, noting that every Hessian matrix $\nabla^2u(x_0)$, with $u: \Om \rga \Rd$ and $x_0 \in \Om$, is in $X(N,d,2)$ when viewed as a bilinear map 

$$(v_1, v_2 ) \mapsto \frac{\partial^2u}{\partial v_1 \partial v_2} (x_0) .$$
We define the cone $\Lambda(N,d,2)$ as

$$\Lambda(N,d,2) := \{a \otimes b \otimes b : a \in \mathbb{R}^d, b \in \mathbb{R}^N \}. $$

\begin{lem}\label{hull}

Let $M = \dim(X(N,d,2)) = \frac{d(d+1)}{2}N$. There is a basis $\{ \xi_i \}_{i=1}^M \subset \Lambda(N,d,2)$ for $X(N,d,2)$ with $|\xi_i| = 1$ for every $i \in \{1, \dots M \}$ and there exists $c(N,d) > 0$ such that for all $H \in X(N,d,2)$ written as
$$H = \sum_{i = 1}^{M} a_i \xi_i, \ a_i \in \mathbb{R}, i = 1, \dots, M, $$
it holds that
$$\frac{1}{c} |H| \leq \sum_{i=1}^{M} |a_i| \leq c |H|. $$

\end{lem}

\proof

Since tensors of the form

$$e_k \otimes e_i  \otimes e_j+ e_k  \otimes e_j  \otimes e_i, \ k = 1, \dots, d, \textrm{ and }  i,j = 1, \dots, N,  $$
\noindent
form a basis for $X(N,d,2)$, to see that we can form a basis contained in $\Lambda(N,d,2)$ it will suffice to show that the span of $\Lambda(N,d,2)$ contains these basis vectors. When $i = j$, we trivially have

$$2 e_k \otimes e_i \otimes e_i \in \Lambda(N,d,2), $$
and if $i \not = j$, we note that $\Lambda(N,d,2)$ contains

$$e_k \otimes (e_i+e_j) \otimes (e_i+e_j) = e_k \otimes ( e_i \otimes e_i +  e_j \otimes e_i +  e_i \otimes e_j +  e_j \otimes e_j  )$$
which, combined with our above observation, implies that

$$e_k \otimes e_i  \otimes e_j+ e_k  \otimes e_j  \otimes e_i \in \textrm{Span}(\Lambda(N,d,2)).$$

Thus we have $\textrm{Span}(\Lambda(N,d,2)) = X(N,d,2)$ and we can select a basis for $X(N,d,2)$ consisting of $\Lambda(N,d,2)$ tensors, and by scaling these appropriately we can guarentee $|\xi_i|= 1$ for every $i \in \{1, \dots, M\}$. Note that with $H = \sum_{i=1}^M a_i \xi_i$, $a_i \in \mathbb{R}$, $i = 1, \dots, M$,
$$\| H \| := \sum_{i=1}^M |a_i| $$
defines a norm on $X(N,d,2)$, and the existence of a constant $c>0$ such that
$$\frac{1}{c} |H| \leq \sum_{i=1}^{M} |a_i| \leq c |H| $$
follows from the equivalence of norms on finite dimensional normed spaces. 

\qed

\bd{Lambda Convex}

We say that a function $F : X(N,d,2) \rightarrow \mathbb{R}$ is $\Lambda(N,d,2)$-convex if
\noindent
$$F(t \xi + (1-t) \xi') \leq t F(\xi) + (1-t) F(\xi') $$
whenever $(\xi - \xi') \in \Lambda(N,d,2)$, $t \in (0,1)$. 

\ed

Theorem 3.3 in \cite{BCO81} relates 2-quasiconvexity to $\Lambda(N,d,2)$-convexity.

\bl{Null Lagrangians}

Let $F : \mathbb{R}^{N \times d \times d} \rga \mathbb{R}$ be continuous and 2-quasiconvex. Then $F$ is $\Lambda(N,d,2)$-convex. 

\el

Next we show that $\Lambda(N,d,2)$-convex functions with linear growth are in fact Lipschitz continuous in all of $X(N,d,2)$. This lemma is a slight modification of a similar result on separately convex functions in \cite{FL07} Proposition 4.64. In this proof we use a mollification argument. Here, and in what follows, $\phi \in C^{\infty}(\Rn; [0,\infty))$ is such that supp$(\phi) \subset B(0,1)$ and $\int_{\mathbb{R}^N} \phi(x) dx = 1$. We define the standard mollifiers $\phi_{\eps}(x) := \frac{1}{\eps^N} \phi(\eps x)$, $\eps > 0$.

\bl{LipCon}

Let $f : \mathbb{R}^{d \times N \times N} \rga \mathbb{R}$ be a $\Lambda(N,d,2)$-convex function such that
\be{lin}
|f(H)| \leq C(1+|H|)
\ee
for some $C > 0$ and all $H \in \mathbb{R}^{d \times N \times N}$. Then
$$|f(H) - f(H')| \leq \tilde{C}|H - H'| $$
for all $H, H' \in X(N,d,2)$, where $\tilde{C}$ depends only on $C,N$ and $d$.
\el

\proof 

\underline{Step 1}: First, consider the case where $f \in C^{\infty}(\mathbb{R}^{N \times d \times d})$. From Lemma \ref{hull} we can select a basis $\{ \xi_i \} \subset \Lambda(N,d,2)$ for $X(N,d,2)$. Fix $H \in X(N,d,2)$, which can be expressed as $H = \sum_{i = 1}^{M} a_i \xi_i, a_i \in \mathbb{R}, i = 1, \dots, M$. Fix $j \in \{1, \dots, M\}$ and consider the function
$$g(t) = f\bigg( t \xi_j +\sum_{i \not = j}a_i \xi_i\bigg).  $$
Since $g$ is convex and smooth, it follows from \cite{FL07} Theorem 4.62 that for every $t, s \in \mathbb{R}$ we have
$$g(t + s) - g(t) \geq g'(t) s.$$
In particular, letting $s = 1 + |H|$ and $t = a_j$,
\begin{align*}
g'(t) = \frac{\partial f}{\partial \xi_j}(H) &\leq \frac{g(t+s)-g(t)}{s} \leq \frac{|f(H+ (1+|H|) \xi_j)| + |f(H)|}{1+|H|}  \\
&= \frac{C(1+|H| + |\xi_j|(1 + |H|)) + C(1+|H|)}{1+|H|} \leq 3C \frac{1 + |H|}{1 + |H|} =3C,
\end{align*}
by virtue of (\ref{lin}) and the fact that $|\xi_j|=1$. Similarly,
$$g(t-s) - g(t) \geq -g'(t)s $$
so
$$-g'(t) \leq \frac{g(t-s)- g(t)}{s} \leq 3C \frac{1+ |H|}{1+|H|} = 3C $$
and thus
$$\bigg| \frac{\partial f}{\partial \xi_j}(H) \bigg| \leq 3C $$
for every $j = 1, \dots, M,$ and $H \in X(N,d,2)$. Let $H, H' \in X(N,d,2)$. By the mean value theorem, we can find $\theta \in (0,1)$ so that
\be{Liph}
|f(H) - f(H')| = |\nabla f ( \theta H + (1- \theta) H') \cdot (H-H')| 
\ee
and we can decompose $H - H'$ into $\sum_{i =1}^{M} b_i \xi_i$ so that
\begin{align*}
|\nabla f ( \theta H + (1- \theta) H') \cdot (H-H')| &= \bigg| \sum_{i = 1}^{M} \sum_{j = 1}^{M} \frac{\partial f}{\partial \xi_i}(\theta H + (1- \theta) H') b_j \xi_i \cdot \xi_j \bigg| \\
& \leq \sum_{i = 1}^M \sum_{j = 1}^M 3C |b_j|  \leq 3cMC |H-H'|
\end{align*}
where we used Lemma \ref{hull}. We conclude in view of (\ref{Liph}).

\underline{Step 2}: For an arbitrary $\Lambda(N,d,2)$-convex function $f$ satisfying (\ref{lin}), consider the mollified functions $f_{\eps} := f * \phi_{\eps}$, $\eps >0$. Each function $f_{\eps}$ is still $\Lambda(N,d,2)$-convex and for every $H \in \mathbb{R}^{N \times d \times d}$ we have
\begin{align*}
|f_{\eps}(H)|  &\leq \bigg|\int_{\mathbb{R}^{N \times d \times d}} \phi_{\eps}(S) f(H-S) dS \bigg| \leq  C \int_{B(0, \eps)} \phi_{\eps}(S) (1+ |H -S|) dS \\
&\leq C (1 + |H|) ,
\end{align*}
and by Step 1
$$|f_{\eps}(H) - f_{\eps}(H')| \leq \tilde{C}|H-H'| $$
for every $H, H' \in X(N,d,2)$ for some $\tilde{C}$ independent of $\eps$. Since $f_{\eps} \rga f$ pointwise as $\eps \rga 0^+$, we have our desired result.

\qed

To prove the upper bound, we will establish an area-strict density result in $BH$. The notion of area-strict convergence, as discussed in \cite{RS15}, is as follows.

\bd{Area Strict Convergence}
We say that a sequence of Radon measures $\{ \mu^n \} \subset \cM(\Om; \mathbb{R}^k)$ \textnormal{converges area-strictly} to $\mu \in \cM(\Om; \mathbb{R}^k)$ if $\mu^n \wlims \mu$, i.e.,
\[\int_{\Om} \psi \cdot d \mu_n \rga \int_{\Om} \psi \cdot d \mu \textrm{ for every } \psi \in C_c(\Om; \mathbb{R}^k),\] 
and
$$\int_{\Om} \sqrt{ 1+  \big|\mu^n_{ac}\big|^2} dx + |\mu^n_s|(\Om) \rightarrow \int_{\Om} \sqrt{ 1+  \big|\mu_{ac}\big|^2} dx + |\mu_s|(\Om).$$
\ed

We will make use of another Reshetnyak-type theorem found in \cite{KR10} Theorem 5.

\begin{thm}\label{eContinuity}
Let $f \in \mathbf{E}(\Om; \mathbb{R}^{d \times N}) $ and let $\{\mu^n\}$ be a sequence of matrix valued measures on $\Om$ such that $\mu^n \rightarrow \mu$ area-strictly on $\Om$. Then

$$\int_{\Om} f(x,  \mu^n_{ac}) dx + \int_{\Om} f^{\infty} \Bigg( x, \frac{d \mu^n_s}{d|\mu^n_s|} \Bigg) d|\mu^n_s| \rightarrow \int_{\Om} f (x, \mu_{ac}) dx + \int_{\Om} f^{\infty} \Bigg(x, \frac{d\mu_s}{d|\mu_s|} \Bigg) d|\mu_s|,$$

\end{thm}
\noindent
where
\begin{align*}
\mathbf{E}(\Om; \mathbb{R}^{d \times N}) := \Big\{ f : \ &\Om \times \mathbb{R}^{d \times N} \rightarrow \mathbb{R} \ | \ \hat{f}(x, \xi) := (1-|\xi|) f(x, (1-|\xi|)^{-1} \xi) 
 \textrm{ as a function} \\
&\textrm{in } \Om \times B(0,1) \textrm{ has a continuous extension to  } \overline{\Om \times B(0,1)} \Big\}.
\end{align*}
We apply Theorem \ref{eContinuity} to obtain the following result:

\begin{thm}\label{Area-Strict Continuity}
Let $f : \Om \times \mathbb{R}^{d \times N \times N} \rga [0, \infty)$ be a 2-quasiconvex continuous integrand satisfying the growth condition (H1). Then the functional

$$\mathcal{G}[u] := \int_{\Om} f(x, \nabla^2u(x)) dx + \int_{\Om} f^{\infty} \Bigg(x, \frac{d D_s(\nabla u)}{d|D_s(\nabla u)|}(x) \Bigg) d|D_s(\nabla u)|(x) $$
is continuous with respect to area-strict convergence of $D(\nabla u)$.

\end{thm}

\proof

From Lemma 1 in \cite{KR10}, we know that since $f$ is continuous and nonnegative with linear growth, we can find $g_k, h_k \in \mathbf{E}(\Om; \mathbb{R}^{d \times N \times N})$ such that

$$g_k(x, H) \nearrow f(x, H), \  g_k^{\infty}(x,H) \nearrow f_{\#}(x, H),$$
$$h_k(x,H) \searrow f(x,H), \  h_k^{\infty}(x,H) \searrow f^{\#}(x,H), $$
for every $x \in \Om$, $H \in X(N,d,2)$, where

$$ f_{\#}(x, H) := \liminf \bigg\{ \frac{f(x',tH')}{t} : x' \rga x, H' \rga H, t \rga + \infty \bigg\},$$
and
$$  f^{\#}(x, H) := \limsup \bigg\{ \frac{f(x',tH')}{t} : x' \rga x, H' \rga H, t \rga + \infty \bigg\}. $$

Let $u \in BH(\Om; \Rd)$ and let $u_n \in W^{2,1}(\Om; \Rd)$ be such that $u_n \rga u$ in $W^{1,1}$ and $\nabla^2 u_n \Ln \mres \Om \rga D(\nabla u)$ area-strictly. For every $k$ we apply Theorem \ref{eContinuity} to obtain

\begin{align*}
\liminf_{n \rga \infty} \cG[u_n] &\geq \liminf_{n \rga \infty} \int_{\Om} g_k(x, \nabla^2 u_n) dx \\
&=\int_{\Om} g_k(x, \nabla^2 u) dx + \int_{\Om} g_k^{\infty}\bigg(x, \frac{d D_s(\nabla u)}{d |D_s(\nabla u)|} \bigg) d|D_s(\nabla u),\
\end{align*}
and

\begin{align*}
\limsup_{n \rga \infty} \cG[u_n] &\leq \limsup_{n \rga \infty} \int_{\Om} h_k(x, \nabla^2 u_n) dx \\
&=\int_{\Om} h_k(x, \nabla^2 u) dx + \int_{\Om} h_k^{\infty}\bigg(x, \frac{d D_s(\nabla u)}{d |D_s(\nabla u)|} \bigg) d|D_s(\nabla u).
\end{align*}

Taking the supremum over $k$, we apply Monotone Convergence to conclude
\be{gs}
 \liminf_{n \rga \infty} \cG[u_n] \geq \int_{\Om} f(x, \nabla^2 u) dx + \int_{\Om} f_{\#}\bigg(x, \frac{d D_s(\nabla u)}{d |D_s(\nabla u)|} \bigg) d|D_s(\nabla u)|.
\ee
Similarly, since $h_1, h_1^{\infty} \in \mathbb{E}(\Om; \mathbb{R}^{d \times N \times N})$ we can apply Monotone Convergence to $-h_k, -h_k^{\infty}$ to conclude
\be{hs}
\limsup_{n \rga \infty} \cG[u_n] \leq  \int_{\Om} f(x, \nabla^2 u) dx + \int_{\Om} f^{\#}\bigg(x, \frac{d D_s(\nabla u)}{d |D_s(\nabla u)|} \bigg) d|D_s(\nabla u)|.
\ee
A generalization of Alberti's Rank One theorem to Hessians, proved in \cite{dPR16}, Theorem 1.6, says that 
\be{albertiR1}
\frac{D_s(\nabla u)}{|D_s(\nabla u)|}(x) \in \Lambda(N,d,2)
\ee
for $|D_s(\nabla u)|$ almost every $x$. We claim that for all $H \in \Lambda(N,d,2)$
\be{sameAlb}
f_{\#}(x, H) = f^{\#}(x, H) = f^{\infty}(x,H) .
\ee
To see this, as in \cite{KR10}, we examine the expression

\begin{align}
\label{H prime}
\frac{f(x', tH')}{t} &= \frac{f(x', tH') - f(x', tH)}{t} + \frac{f(x', 0)}{t} + \frac{f(x', tH) - f(x', 0)}{t}
\end{align}
for $x' \in \Om$, $H' \in X(N,d,2)$ and $t > 0$. By Lemma \ref{Null Lagrangians}, $f(x, \cdot)$ is $\Lambda(N,d,2)$-convex with linear growth. Hence by Lemma \ref{LipCon}, $f(x, \cdot)$ is Lipschitz on all of $X(N,d,2)$ and the Lipschitz constant is independant of $x$. Thus, the first term in (\ref{H prime}) will vanish as $H' \rga H$.

The second term clearly goes to zero as $t \rga + \infty$, so we turn to the third term. We note that for all $H \in \Lambda(N,d,2)$, $y \in \Om$, 

$$ \frac{f(y, tH) - f(y,0)}{t}$$
is an increasing function in $t$ by $\Lambda(N,d,2)$-convexity, and since $f(y, \cdot)$ is Lipschitz, we have, recalling (\ref{recdef}),
$$\lim_{t \rga + \infty} \frac{f(y, tH) - f(y,0)}{t} = \sup_{t > 0} \frac{f(y, tH) - f(y,0)}{t} = f^{\infty}(y, H).$$
As $f$ is continuous in $y$ for every $H$, we can apply Dini's Theorem to conclude that the convergence as $t \rga + \infty$ is locally uniform in $y$. Thus, the third term converges to $ f^{\infty}(x, H)$ as $t \rga + \infty$ and $x' \rga x$.

In view of (\ref{gs}), (\ref{hs}), (\ref{albertiR1}) and (\ref{sameAlb}) we conclude that

$$\lim_{n \rga \infty} \cG[u_n] = \cG[u]. $$

\qed

\section{Density Result}

Here we prove a useful density result which states that we can approximate a measure in the area-strict sense via smooth functions, as long as the domain is sufficiently regular. In order to prove this, we will need the following estimates.

\bl{convolutions}

Let $\Om \subset \Rn$ be an open set. Let $g: \Rd \rga \mathbb{R}$ be a convex function satisfying $|g(\xi)| \leq C(1 + |\xi|)$ for some $C>0$ and all $\xi \in \Rd$, and let $\mu \in \cM(\Om, \Rd)$. For every $x \in \Omega$ and $\eps < \eps_0 :=  \textnormal{dist}(x, \partial \Om)$,

$$
g((\mu_{ac} * \phi_{\eps} )(x))  \leq (g(\mu_{ac}) * \phi_{\eps})(x),
$$
and
\begin{equation}
\label{singBound}
g(( \mu_{s} * 2 \phi_{\eps}) (x) ) \leq \int_{\Omega} \frac{g(2 t_{\eps}(x) \nu(y))}{t_{\eps}(x)} \phi_{\eps}(y-x) d|\mu_s|(y),
\end{equation}
where $t_{\eps} \in C^{\infty}(B(x, \eps)) ; [0, \infty))$ is given by
$$t_{\eps}(x) = \int_{\Omega} \phi_{\eps} (y-x) d |\mu_{s}|(y) dy,$$
and (\ref{singBound}) holds whenever $t_{\eps}(x)>0$, a set of $|\mu_s|$ density 1.

\el

\proof

Fix $x \in \Om$ and note that $\eps < \eps_0$ implies that $B(x, \eps) \subset \Omega$. By Jensen's inequality,

$$ g( (\mu_{ac} * \phi_{\eps} )(x)) = g \bigg( \int_{\Omega} \phi_{\eps} (y-x) \mu_{ac}(y) dy \bigg) \leq \int_{\Omega} \phi_{\eps} (y-x) g(\mu_{ac}(y)) dy = (g(\mu_{ac}) * \phi_{\eps}) (x).$$
\noindent
where we used the fact that $\int_{\Omega} \phi_{\eps}(y-x) dy = 1$.

\

For the singular part, we set $t_{\eps}(x) := \int_{\Omega} \phi_{\eps} (y-x) d |\mu_{s}|(y) dy$. If $\eps <<1$, then $t_{\eps} \in C^{\infty}(B(x, \eps_0); [0, \infty))$, and if $t_{\eps} > 0$ then the measure $\pi_{\eps} := \frac{1}{t_{\eps}} \phi_{\eps}(\cdot - x) |\mu_s|$ is a probability measure. Thus, we can again apply Jensen's inequality to obtain

\begin{align*}
g( (\mu_{s} * 2 \phi_{\eps})(x)) &= g \bigg( \int_{\Omega} 2 \phi_{\eps} (y-x) d \mu_{s}(y)  \bigg) = g \bigg( \int_{\Omega} 2 t_{\eps}(x) \nu(y) d\pi_{\eps}(y) \bigg)\\
&\leq \int_{\Omega} g(2 t_{\eps}(x) \nu(y)) d \pi_{\eps}(y) = \int_{\Omega} \frac{g(2 t_{\eps}(x) \nu(y))}{t_{\eps}(x)} \phi_{\eps}(y-x) d|\mu_s|(y) .
\end{align*}

\qed

In order to establish the lower bound for the singular part in both the density result and the relaxation result, we invoke a lemma which can be found in \cite{FM93}, Lemma 2.13. To be precise,

\bl{smallercubes}

Let $\lambda$ be a nonnegative Radon measure on $\Rn$. For $\lambda$ almost every $x_0 \in \Rn$ and for every $0 < \sigma <1$, 

\begin{equation}
\label{FMcubes}
\limsup_{r \rga 0^+} \frac{\lambda(Q(x_0, \sigma r))}{\lambda(Q(x_0, r))} \geq \sigma^N.
\end{equation}

\el 
We will use a modification of this lemma: in (\ref{FMcubes}) we can choose $r\rga 0^+$ so that, given another Radon measure $\mu$, neither $\mu$ nor $\lambda$ charge the boundary of the larger cubes. Namely, we have the following result.

\bl{smallercubes2}

Let $\lambda$ and $\mu$ be nonnegative Radon measures in $\Rn$. For every $0 < \sigma < 1$, and for $\lambda$ almost every $x_0 \in \Rn$, there exist $r_n \rga 0^+$ such that $\mu(\partial Q(x_0, r_n)) = \lambda(\partial Q(x_0, r_n)) = 0$ and

$$\lim_{n \rga \infty} \frac{\lambda(Q(x_0, \sigma r_n))}{\lambda(Q(x_0, r_n))} \geq \sigma^N .$$
 
\el

\proof

Fix $\sigma \in (0,1)$ and $x_0 \in \Rn$ so that, by Lemma \ref{smallercubes}, we can find $\rho_n \rga 0^+$ such that

\be{goodR}
\lim_{n \rga \infty} \frac{\lambda(Q(x_0, \sigma \rho_n))}{\lambda(Q(x_0, \rho_n))} \geq \sigma^N.
\ee
For every $n \in \mathbb{N}$, we can select $\delta_n < \rho_n$ such that

$$\lambda(Q(x_0, \sigma \delta_n)) \geq \frac{n}{n+1} \lambda(Q(x_0, \sigma \rho_n)). $$
Find $\tilde{r}_n \in (\delta_n, \rho_n)$ such that $\mu(\partial Q(x_0, \tilde{r}_n)) = \lambda(\partial Q(x_0, \tilde{r}_n)) = 0$. We obtain

$$ \frac{\lambda(Q(x_0, \sigma \tilde{r}_n))}{\lambda(Q(x_0, \tilde{r}_n))} \geq \frac{\lambda(Q(x_0, \sigma \delta_n))}{\lambda(Q(x_0, \rho_n))} \geq \frac{n}{n+1} \frac{\lambda(Q(x_0, \sigma \rho_n))}{\lambda(Q(x_0, \rho_n))},$$
and by (\ref{goodR}) we conclude that

$$\liminf_{n \rga \infty} \frac{\lambda(Q(x_0, \sigma \tilde{r}_n))}{\lambda(Q(x_0, \tilde{r}_n))} \geq \sigma^N. $$
Since the sequence $\big\{ \frac{\lambda(Q(x_0, \sigma \tilde{r}_n))}{\lambda(Q(x_0, \tilde{r}_n))} \big\} $ is bounded, we can extract a subsequence $\{r_n\}_{n \in \mathbb{N}} \subset \{\tilde{r}_n\}_{n \in \mathbb{N}}$ such that

$$\lim_{n \rga \infty} \frac{\lambda(Q(x_0, \sigma r_n))}{\lambda(Q(x_0, r_n))} = \liminf_{n \rga \infty} \frac{\lambda(Q(x_0, \sigma \tilde{r}_n))}{\lambda(Q(x_0, \tilde{r}_n))} \geq \sigma^N .$$

\qed

We now state and prove the main result of this section.

\bt{density}
\label{density}

Let $U$ be an open, bounded set in $\Rn$, let $\mu \in \cM(U; \Rd)$ be a Radon measure, and let $g: \Rd \rga [0, \infty)$ be a convex function such that
\be{A1}
0 \leq g(p) \leq C (1 + |p|) \textrm{ for all } p \in \mathbb{R}^d \textrm{ and some } C>0
\ee
and
\be{A3}
\bigg|\frac{g(tp)}{t} - g^{\infty}(p)\bigg| \leq \frac{C}{t^{\alpha}}
\ee
for some $\alpha >1$ and $C >0$ and all $p \in \Rd$ with $|p|=1$ and every $t >0$. Then, for every $\Om \subset \subset U$ with $|\partial \Omega| = |\mu|( \partial \Omega) = 0$,
$$ 
\lim_{\eps \rga 0^+} \int_{\Omega} g(\mu_{\eps}) dx = \int_{\Omega} g(\mu_{ac}) dx
 + \int_{\Om} g^{\infty}(\nu) d|\mu_s|$$

\noindent
where $\mu_{\eps} := \mu * \phi_{\eps}$ for $\eps >0$.

\et

\proof




\underline{Step 1: Lower bound}. We claim that

\begin{equation}
\int_{\Omega} g(\mu_{ac}) dx
 + \int_{\Om} g^{\infty}(\nu) d|\mu_s| \leq \liminf_{\eps \rga 0^+} \int_{\Omega} g(\mu_{\eps}) dx.
\label{lub}
\end{equation}

For this inequality, we use the fact that $\{ \mu_{\eps} \Ln \mres \Om \}$ converges weakly-$*$ to $\mu$, and that $\{ |\mu_{\eps} \Ln \mres \Om | \}$ converges weakly-$*$ to $|\mu|$ (See \cite{AFP00}, Theorem 2.2). It should be noted that in this step we do not need the assumption that $|\mu|(\partial \Om) = 0$.

We will apply the blow-up argument originally found in \cite{FM93}. Choose $\eps_k \rga 0$ which achieve the liminf, and, for simplicity, using the notation $\mu_{\eps_k} =: \mu_k$, we define the Radon measures

$$\lambda_k(E) :=  \int_{E} g(\mu_k(x)) dx$$
\noindent
for any Borel set $E \subset \Om$. Due to the growth condition (\ref{A1}), we have

\begin{equation}
\lambda_k(\Om) = \int_{\Om} g(\mu_k(x)) dx \leq \int_{\Om} C ( 1 + |\mu_k(x)|) dx = C \bigg( |\Om| + \int_{\Om} |\mu_k| dx \bigg),
\label{Lbound}
\end{equation}
\noindent
and since $\{ \mu_k \Ln \mres \Om \}$ converge weakly-$*$, the sequence $\{ \int_{\Om} | \mu_k | dx \}$ is bounded. We deduce that $\{\lambda_k(\Om)  \}$ is bounded, therefore, along a subsequence (not relabled) we have $\lambda_k \wlims \lambda$ for some nonnegative finite Radon measure $\lambda$. 

\

The growth conditions on $g$ yield 

\begin{equation}
\lambda << \Ln \mres \Om + |\mu|.
\label{lambdaAC}
\end{equation}
Indeed, let $E$ be any Borel subset of $\Om$ with $|E| = |\mu|(E) = 0$. By inner regularity, it suffices to show $\lambda(K) = 0$ for every $K \subset E$ compact. For any such $K$, we have 
\begin{equation}
|K| = |\mu|(K) =0.
\label{CompactAC}
\end{equation}
Define the open sets
\[ K_{\delta} := \{x \in \Om : \textrm{dist}(x, K) < \delta \}. \]
Since $\partial K_{\delta} = \{ x \in \Om : \textrm{dist}(x,K) = \delta\}$ are an uncountable family of disjoint sets, we can select $\delta_i \rga 0^+$ such that 

\begin{equation}
|\mu|(\partial K_{\delta_i}) = 0.
\label{Kbound}
\end{equation}
We have by ($\ref{Lbound}$)

$$\lambda(K) \leq \lambda(K_{\delta_i}) \leq \liminf_{k \rga \infty} \lambda_k(K_{\delta_i}) \leq \liminf_{k \rga \infty} C \bigg( |K_{\delta_i}| + \int_{ K_{\delta_i}} |\mu_k| dx \bigg) $$

$$= \lim_{k \rga \infty} C \bigg( |K_{\delta_i}| +  \int_{ K_{\delta_i}} |\mu_k| dx \bigg) =  C \bigg( |K_{\delta_i}| + |\mu|(K_{\delta_i}) \bigg) $$
\noindent
by virtue of ($\ref{Kbound}$) and the fact that $|\mu_k| \Ln \mres \Om \wlims |\mu|$. Since $\bigcap_{\delta >0} K_{\delta} = K$, letting $i \rga \infty$, we get by (\ref{CompactAC}) 

$$\lambda(K) \leq C \bigg( |K| + |\mu|(K) \bigg) = 0, $$
and this concludes (\ref{lambdaAC}).

\

We claim that

\begin{equation}
\frac{d \lambda}{d \Ln} (x_0) \geq g(\mu_{ac}(x_0)) \textrm{ for } \Ln a.e. \ x_0 \in \Om,
\label{LBAC}
\end{equation}
and

\begin{equation}
\frac{d \lambda}{d |\mu_s|} (x_0) \geq g^{\infty}(\nu(x_0)) \textrm{ for } |\mu_s|\ a.e. \ x_0 \in \Om.
\label{LBS}
\end{equation}
If ($\ref{LBAC}$) and ($\ref{LBS}$) hold, then
$$\int_{\Om} g(\mu_{ac}(x))dx + \int_{\Om} g^{\infty}(\nu(x)) d |\mu_s|(x)  \leq \lambda(\Om) \leq \liminf_{k \rga \infty}\lambda_k( \Om) = \lim_{k \rga \infty} \int_{\Om} g( \mu_k(x))dx $$
and this yields (\ref{lub}).

We begin by establishing the inequality for the absolutely continuous part, i.e., (\ref{LBAC}). For $\Ln$ almost every $x_0 \in \Om$, we have
$$\frac{d \lambda}{d \Ln} (x_0) = \lim_{r \rga 0} \frac{\lambda(Q(x_0, r))}{r^N} \textrm{ exists and is finite,}$$

\begin{equation}
\lim_{r \rga 0} \frac{1}{r^N} \int_{Q(x_0, r)} |\mu_{ac}(x) - \mu_{ac}(x_0)| dx = 0,
\label{A}
\end{equation}

\begin{equation}
\lim_{r \rga 0} \frac{|\mu|(\overline{Q(x_0,r)})}{r^N} = |\mu|_{ac}(x_0) = |\mu_{ac}|(x_0),
\label{C} 
\end{equation}
and
\begin{equation}
\lim_{r \rga 0} \frac{|\mu_s|(\overline{Q(x_0, r)})}{r^N} = \lim_{r \rga 0} \frac{|\mu|_s(\overline{Q(x_0, r)})}{r^N} = 0.
\label{B}
\end{equation}
where we have used Lemma \ref{Radon Nikodym} in (\ref{C}) and (\ref{B}). Choose $r_n \rga 0$ such that $|\mu|(\partial Q(x_0, r_n)) = \lambda(\partial Q(x_0, r_n)) = 0$. Then 

$$\frac{d \lambda}{d \Ln}(x_0) = \lim_{n \rga \infty} \frac{\lambda(Q(x_0, r_n))}{r_n^N} = \lim_{n \rga \infty} \lim_{k \rga \infty} \frac{\lambda_k(Q(x_0, r_n))}{r_n^N}$$

\begin{equation}
= \lim_{n \rga \infty} \lim_{k \rga \infty} \frac{1}{r_n^N} \int_{Q(x_0, r_n)} g(\mu_k(x)) dx .
\label{FMKAC}
\end{equation}

Define the functions $v_{n,k}(y) := \mu_k(x_0 + r_n y)$ for $y \in Q.$ Apply a change of variables so that ($\ref{FMKAC}$) becomes
\be{CoV}
\frac{d \lambda}{d \Ln}(x_0) = \lim_{n \rga \infty} \lim_{k \rga \infty} \int_{Q} g( v_{n, k}(y)) dy. 
\ee
Since $\mu_k \Ln \mres \Om \wlims \mu$, we have that for every $n \in \mathbb{N}$, the measures $v_{n,k} \ \Ln \mres Q$ converge weakly-$*$ to a measure $\pi_n$ given by

\begin{equation}
\pi_n(E) := \frac{(T^{\#}_{x_0, r_n} \mu) (E)}{r_n^N}  = \frac{\mu(x_0 + r_n E)}{r_n^N}, \textrm{ for every Borel set } E \subset Q, 
\label{PFAC}
\end{equation}
where $T^{\#}_{x_0, r_n} \mu$ denotes the push-forward of $\mu$ under the mapping $x \mapsto \frac{x-x_0}{r_n}.$ Indeed, by the standard change of variables for push-forward measures (see \cite{B07}, Theorem 3.6.1), for any test function $\psi \in C_c(Q)$ we have
\begin{align}
\notag
\int_{Q} v_{n, k}(y) \psi(y) dy &= \int_{Q} \mu_k(x_0 + r_n y) \psi(y) dy = \frac{1}{r_n^N} \int_{Q(x_0, r_n)} \mu_k(x) \psi \bigg( \frac{x-x_0}{r_n} \bigg) dx\\ 
\notag
&\buildrel{k \rga 0}\over\longrightarrow  \ \frac{1}{r_n^N} \int_{Q(x_0, r_n)} \psi \bigg( \frac{x-x_0}{r_n} \bigg) d\mu(x) = \frac{1}{r_n^N} \int_{Q} \psi(y) d (T^{\#}_{x_0, r_n} \mu) (y) \\
&= \int_Q \psi(y) d \pi_n(y).
\notag
\end{align}
In turn, $\pi_n \wlims \mu_{ac}(x_0) \Ln \mres Q$. To see this, fix any $\psi \in C_c(Q)$. We have

$$\bigg|\int_{Q} \psi(y) d \pi_n(y) - \int_{Q} \psi(y) \mu_{ac}(x_0)  dy \bigg|$$
\begin{align*}
&= \frac{1}{r_n^N} \bigg|\int_{Q(x_0, r_n)} \psi\bigg(\frac{x-x_0}{r_n}\bigg) d\mu(x) - \int_{Q(x_0, r_n)} \psi\bigg(\frac{x-x_0}{r_n}\bigg) \mu_{ac}(x_0) dx \bigg|\\
&= \frac{1}{r_n^N} \bigg| \int_{Q(x_0, r_n)} \psi\bigg(\frac{x-x_0}{r_n}\bigg) (\mu_{ac}(x)-\mu_{ac}(x_0)) dx + \int_{Q(x_0, r_n)} \psi\bigg(\frac{x-x_0}{r_n}\bigg) d\mu_s(x) \bigg|\\
&\leq \|\psi\|_{\infty} \bigg( \frac{1}{r_n^N} \int_{Q(x_0, r_n)} |\mu_{ac}(x) - \mu_{ac}(x_0)| dx + \frac{|\mu_s|(Q(x_0, r_n))}{r_n^N} \bigg)
\end{align*}
which goes to $0$ as $n \rga \infty$ by $(\ref{A})$ and $(\ref{B})$. Thus, $\pi_n \wlims \mu_{ac}(x_0) \ \Ln \mres Q$. We also observe
\begin{align}
\notag
\limsup_{n \rga \infty} \limsup_{k \rga \infty} \|v_{n,k} \Ln \mres Q \| &= \limsup_{n \rga \infty} \limsup_{k \rga \infty} \int_{Q} |\mu_k|(x_0 + r_n y) dy \\
\notag
&= \limsup_{n \rga \infty} \frac{1}{r_n^N} \limsup_{k \rga \infty} \int_{Q(x_0, r)} |\mu_k|(x) dx \\
\notag
&\leq \limsup_{n \rga \infty} \frac{|\mu|(\overline{Q(x_0,r)})}{r^N} = |\mu|_{ac}(x_0).
\end{align}

If we fix $t \in (0,1)$ and note that

\begin{align*}
\notag
\varlimsup_{n \rga \infty} \varlimsup_{k \rga \infty} \int_{Q \setminus tQ} |v_{n,k}(y)| dy &= \varlimsup_{n \rga \infty} \varlimsup_{k \rga \infty} \frac{1}{r_n^N} \int_{Q(x_0, r_n) \setminus Q(x_0, t r_n)} |\mu_k(x)|dx \\ 
\notag
&\leq \varlimsup_{n \rga \infty} \frac{|\mu|(Q(x_0, r_n) \setminus Q(x_0, t r_n))}{r_n^N} \\
\notag
&= \lim_{n \rga \infty} \bigg( \frac{|\mu|(Q(x_0, r_n))}{r_n^N} - t^N \frac{|\mu|(Q(x_0, t r_n))}{t^N r_n^N} \bigg) \\
\notag
& = |\mu|_{ac}(x_0)(1-t^N),
\end{align*}
then by Lemma \ref{diags} we can find a diagonal sequence $v_n := v_{n, k_n}$ such that  
\be{diag1}
\lim_{n \rga \infty} \lim_{k \rga \infty} \int_{Q} g( v_{n, k}(y)) dy  = \lim_{n \rga \infty} \int_{Q} g( v_n(y)) dy,
\ee
\be{vgood}
\varlimsup_{n \rga \infty} \int_{Q \setminus tQ} |v_n(y)|dy \leq |\mu_{ac}(x_0)|(1-t^N),
\ee
and $v_n \wlims \mu_{ac}(x_0) \ \Ln \mres Q$.

Since $g$ is convex, in view of Theorem 5.14 in \cite{FL07} consider an affine function $a + b \cdot \xi \leq g(\xi)$ and observe that by (\ref{CoV}) and (\ref{diag1}),
\begin{align*}
\notag
\frac{d \lambda}{d \Ln} (x_0) &= \lim_{n \rga \infty} \int_Q g(v_n(y))dy \geq \liminf_{n \rga \infty} \int_Q \bigg( a + b \cdot v_n(y) \bigg) dy \\
&= a + \liminf_{n \rga \infty} \int_Q b \cdot v_n(y) dy.
\end{align*}
Let $\psi_t \in C_c(Q; [0,1])$ be such that $\psi_t = 1$ in $tQ$. We have
\begin{align*}
\liminf_{n \rga \infty} \int_Q b \cdot v_n(y) dy &= \liminf_{n \rga \infty} \bigg( \int_Q b \cdot v_n \psi_t dy + \int_Q b \cdot v_n (1-\psi_t) dy \bigg) \\
& \geq \int_Q b \cdot \mu_{ac}(x_0) \psi_t  dy - \limsup_{n \rga \infty} |b| \int_{Q \setminus tQ} |v_n| dy \\
& \geq b \cdot \mu_{ac}(x_0) t^N - |b| |\mu_{ac}(x_0)|(1-t^N),
\end{align*}
where we have used (\ref{vgood}). Thus for any $t < 1$ we have
$$\frac{d \lambda}{d \Ln} (x_0) \geq a + b \cdot \mu_{ac}(x_0) t^N - |b| |\mu_{ac}(x_0)|(1-t^N),$$
and sending $t \rga 1^-$ yields
$$\frac{d \lambda}{d \Ln} (x_0) \geq a + b \cdot \mu_{ac}(x_0).$$
Since this holds for any affine function below $g$, we conclude that

$$\frac{d \lambda}{d \Ln} (x_0) \geq g( \mu_{ac}(x_0)). $$




To address the singular part, we fix $\sigma \in (0,1)$. The measures $|\mu_k| \Ln \mres \Om$ are bounded in total variation, so along a subsequence, not relabeled, they converge weakly-$*$ to some nonnegative Radon measure $\tau$. We can decompose $\tau$ into measures $\tau^A$ and $\tau^B$ such that $\tau^A << |\mu_s|$ and $\tau^B \perp |\mu_s|$. Now, for $|\mu_s|$ almost every $x_0$,

$$ \frac{d \lambda}{d |\mu_s|} (x_0) = \lim_{r \rga 0} \frac{\lambda(\overline{Q(x_0, r)})}{|\mu_s|(\overline{Q(x_0, r)})} = \lim_{r \rga 0} \frac{\lambda(\overline{Q(x_0, r)})}{|\mu|(\overline{Q(x_0, r)})}, $$

\begin{equation}
\lim_{r \rga \infty} \frac{1}{|\mu_s|(\overline{Q(x_0, r)})} \int_{Q(x_0, r)} |\mu_{ac}|(x) dx= 0,
\label{acsmaller}
\end{equation}

\begin{equation}
\intav_{\overline{Q(x_0, r)}} | \nu(x) - \nu(x_0) | d|\mu_s|(x) = 0,
\label{nufact}
\end{equation}

\be{tauA}
\frac{d \tau^A}{d |\mu_s|}(x_0) = \lim_{r \rga 0} \frac{\tau(\overline{Q(x_0, r)})}{|\mu_s|(\overline{Q(x_0, r)})} = \lim_{r \rga 0} \frac{\tau(\overline{Q(x_0, r)})}{|\mu|(\overline{Q(x_0, r)})} < \infty,
\ee

and by Lemma \ref{smallercubes2} we may select $r_n \rga 0$ such that $|\mu|( \partial Q(x_0, r_n))= \tau(\partial Q(x_0, r_n))= 0$ and

\begin{equation}
\lim_{n \rga \infty} \frac{|\mu|(Q(x_0, \sigma r_n))}{|\mu|(Q(x_0, r_n))} \geq \sigma^N.
\label{mostcube}
\end{equation}

Note that in view of (\ref{lambdaAC}), $\lambda(\partial Q(x_0, r_n)) = 0$ for all $n \in \mathbb{N}$. We have 

\begin{align}
\notag
\frac{d \lambda}{d |\mu_s|}(x_0) &= \lim_{n \rga \infty} \frac{\lambda(Q(x_0, r_n))}{|\mu|(Q(x_0, r_n))} = \lim_{n \rga \infty} \lim_{k \rga \infty} \frac{\lambda_k(Q(x_0, r_n))}{|\mu|(Q(x_0, r_n))}\\ 
\label{fmks}
&= \lim_{n \rga \infty} \lim_{k \rga \infty} \frac{1}{|\mu|(Q(x_0, r_n))} \int_{Q(x_0, r_n)} g( \mu_k(x)) dx.
\end{align}
Let 

$$ t_n := \frac{|\mu|(Q(x_0, r_n))}{r_n^N},$$ 
and define 

$$v_{n,k}(y) := \frac{\mu_k(x_0 + r_n y)}{|\mu|(Q(x_0, r_n))} r_n^N = \frac{\mu_k(x_0 + r_n y)}{t_n}.$$
We apply a change of variables to ($\ref{fmks}$) to get
\[ \frac{d \lambda}{d |\mu_s|}(x_0) = \lim_{n \rga \infty} \lim_{k \rga \infty} \int_Q \frac{1}{t_n} g( t_n v_{n,k}(y)) dy. \]
For a fixed $n$,  as $k \rga \infty$ we have $v_{n,k} \Ln \mres Q \wlims \pi_n,$ with

\[ \pi_n(E) := \frac{(T^{\#}_{x_0, r_n} \mu) (E)}{|\mu|(Q(x_0, r_n))}  = \frac{\mu(x_0 + r_n E)}{|\mu|(Q(x_0, r_n))} \textrm{, for every Borel set } E \subset Q, \]
where, as in (\ref{PFAC}), $T^{\#}_{x_0, r_n} \mu$ denotes the push-forward of $\mu$ under the mapping $x \rga \frac{x-x_0}{r_n}.$ We note also that
$$|v_{n,k}| = \frac{|\mu_k|(x_0 + r_n y)}{|\mu|(Q(x_0, r_n))} r_n^N, $$
and hence $|v_{n,k}|\Ln \mres Q \wlims \tau_n$ where
$$\tau_n(E) := \frac{\tau(x_0+r_nE)}{|\mu|(Q(x_0, r_n))}. $$
Define the measures
$$\rho_n(E) := \frac{|\mu|(x_0 + r_n E)}{|\mu|(Q(x_0, r_n))} \textrm{, for every Borel set } E \subset Q .$$
Then $\pi_n \wlims \pi$, $\rho_n \wlims \rho$, for a Radon measure $\pi \in \cM(Q; \Rd)$ and $\rho$ a finite nonnegative Radon measure in $\cM(Q)$, perhaps along a subsequence (not relabeled). We claim that 

\be{normal}
\pi = \nu(x_0) \rho.
\ee
Indeed, fix $\psi \in C_c(Q)$. By Lemma \ref{Radon Nikodym} we have

$$ \bigg| \int_{Q} \psi(y) d \pi_n(y) - \int_{Q} \psi(y) \nu(x_0) d \rho_n(y) \bigg| $$
\begin{align*}
&=\frac{1}{|\mu|(Q(x_0, r_n))} \bigg| \int_{Q(x_0, r_n)} \psi\bigg(\frac{x-x_0}{r_n}\bigg) d \mu(x) - \int_{Q(x_0, r_n)} \psi\bigg(\frac{x-x_0}{r_n}\bigg) \nu(x_0) d |\mu|(x) \bigg|\\
&\leq \frac{1}{|\mu|(Q(x_0, r_n))} \bigg| \int_{Q(x_0, r_n)} \psi\bigg(\frac{x-x_0}{r_n}\bigg) \nu(x) d |\mu_s|(x) - \int_{Q(x_0, r_n)} \psi\bigg(\frac{x-x_0}{r_n}\bigg) \nu(x_0) d |\mu|_s(x) \bigg|\\
&\ \ \ +\frac{1}{|\mu|(Q(x_0, r_n))} \bigg( \int_{Q(x_0, r_n)} |\psi|\bigg(\frac{x-x_0}{r_n}\bigg) |\mu_{ac}|(x) dx + \int_{Q(x_0, r_n)} |\psi|\bigg(\frac{x-x_0}{r_n}\bigg) |\mu|_{ac}(x) dx \bigg)\\
&\leq \frac{\| \psi \|_{\infty}}{|\mu|(Q(x_0, r_n))}   \int_{Q(x_0, r_n)} |\nu(x) - \nu(x_0)| d|\mu_s|(x) + \frac{2 \| \psi \|_{\infty}}{|\mu|(Q(x_0, r_n))} \int_{Q(x_0, r_n)} |\mu|_{ac}(x) dx,
\end{align*}
which goes to 0 as $n \rga \infty$ in view of (\ref{acsmaller}) and (\ref{nufact}). Since

$$\int_{Q} \psi d\pi_n \rga \int_{Q} \psi d\pi ,$$ 
and

$$\int_{Q} \psi \nu(x_0)  d\rho_n  \rga \int_{Q} \psi  \nu(x_0) d\rho,$$
we conclude that $\pi = \nu(x_0) \rho$.
We note that $\rho(Q) \geq \sigma^N$. To see this, by ($\ref{mostcube}$) we have

\begin{align}
\notag
\rho(Q) \geq \rho(\sigma \overline{Q}) \geq \limsup_{n \rga \infty} \rho_n(\sigma \overline{Q}) &\geq \limsup_{n \rga \infty} \rho_n(\sigma Q) \\ 
&= \limsup_{n \rga \infty} \frac{|\mu|(Q(x_0,\sigma r_n))}{|\mu|(Q(x_0, r_n))} \geq \sigma^N.
\label{sigA}
\end{align}
Now, by (\ref{tauA}),

\begin{align*}
\varlimsup_{n \rga \infty} \varlimsup_{k \rga \infty} \| |v_{n,k}| \Ln \mres Q \| &= \varlimsup_{n \rga \infty} \varlimsup_{k \rga \infty} \int_{Q} \frac{r_n^N}{|\mu|(Q(x_0,r_n))} |\mu_k|(x_0+r_ny) dy \\
&= \varlimsup_{n \rga \infty} \frac{1}{|\mu|(Q(x_0,r_n))} \varlimsup_{k \rga \infty} \int_{Q(x_0, r_n)} |\mu_k|(x) dx \\
&\leq \varlimsup_{n \rga \infty} \frac{\tau(\overline{Q(x_0, r_n)})}{|\mu|(Q(x_0,r_n))} = \frac{d \tau^A}{d |\mu_s|}(x_0) < \infty, 
\end{align*}
hence by Lemma \ref{diags} we may diagonalize as in the absolutely continuous case, setting $v_n := v_{n,k_n}$, obtaining
$$v_n \Ln \mres Q \wlims \nu(x_0) \rho$$
\be{conv1}
\frac{d \lambda}{d |\mu_s|}(x_0) = \lim_{n \rga \infty} \int_Q \frac{1}{t_n} g(t_n v_n(y)) dy,
\ee
where we have used (\ref{mostcube}). Fix $\eta > 0$. For any $p \in \Rd$ and $t >0$, we can apply (\ref{A3}) to $\frac{p}{|p|}$ and $t|p|$ which yields

$$\frac{g(tp)}{t|p|} \geq g^{\infty}\bigg(\frac{p}{|p|}\bigg) - \frac{C}{|p|^{\alpha}t^{\alpha}}$$
and thus
$$\frac{g(tp)}{t} \geq g^{\infty}(p) - \frac{C}{|p|^{\alpha-1} t^{\alpha}} $$
which implies that for $|p|> \Big( \frac{C}{\eta}^{\frac{1}{\alpha-1}} \Big) |t|^{ \frac{-\alpha}{\alpha-1}}$ we have
$$\frac{g(tp)}{t} \geq g^{\infty}(p) - \eta.$$
Define the sets $E_n := \{|v_n| > \Big( \frac{C}{\eta} \Big)^{\frac{1}{\alpha-1}} |t_n|^{ \frac{-\alpha}{\alpha-1}} \}$. Then,

$$\lim_{n \rga \infty} \int_Q \frac{1}{t_n} g(t_n v_n(y)) dy \geq \limsup_{n \rga \infty} \int_{E_n} \frac{1}{t_n} g(t_n v_n(y)) dy \geq \limsup_{n \rga \infty} \int_{E_n} g^{\infty}(v_n(y)) dy - \eta. $$
On the other hand,

$$\limsup_{n \rga \infty} \int_{E_n} g^{\infty}(v_n(y)) dy \geq \liminf_{n \rga \infty} \int_{Q} g^{\infty}(v_n(y)) dy - \limsup_{n \rga \infty} \int_{Q \setminus E_n} g^{\infty}(v_n(y)) dy $$
and, since

$$ \limsup_{n \rga \infty} \int_{Q \setminus E_n} g^{\infty}(v_n(y)) dy \leq \limsup_{n \rga \infty} \bigg( \frac{C}{\eta} \bigg)^{\frac{1}{\alpha-1}} |t_n|^{ \frac{-\alpha}{\alpha-1}} = 0,$$
because $t_n \rga \infty$ and $\alpha>1$ implies $\frac{-\alpha}{\alpha-1}<-1,$ 
$$\lim_{n \rga \infty} \int_Q \frac{1}{t_n} g(t_n v_n(y)) dy \geq \liminf_{n \rga \infty} \int_Q g^{\infty}(v_n(y)) dy - \eta $$
for every $\eta > 0$, and so
\be{recc}
\lim_{n \rga \infty} \int_Q \frac{1}{t_n} g(t_n v_n(y)) dy \geq \liminf_{n \rga \infty} \int_Q g^{\infty}(v_n(y)) dy.
\ee

Finally, by Theorem \ref{reshet}, since $g^{\infty}$ is convex and 1-homogeneous with the appropriate growth condition, we have lower semicontinuity with respect to weak-$*$ convergence, and so by (\ref{sigA}), (\ref{conv1}) and (\ref{recc}) we obtain

$$\frac{d \lambda}{d |\mu_s|}(x_0) \geq \int_Q g^{\infty}(\nu(x_0)) d \rho(y) \geq g^{\infty}(\nu(x_0)) \sigma^N, $$
and letting $ \sigma \rga 1^{-}$ we conclude that

$$\frac{d \lambda}{d |\mu_s|}(x_0) \geq g^{\infty}(\nu(x_0)).$$

\noindent\makebox[\linewidth]{\rule{\paperwidth}{0.4pt}}

\underline{Step 2: Upper bound}. We claim that

\begin{equation}
\limsup_{\eps \rga \infty} \int_{\Omega} g(\mu_{\eps}) dx \leq \int_{\Omega} g(\mu_{ac}) dx
 + \int_{\Om} g(\nu) d|\mu_s|.
\label{lemUB}
\end{equation}
 
We will use the blow-up method. Choose a sequence $\{\eps_k\}$ which achieves the limsup, and define the measures

$$\lambda_{k} (E) := \int_{E} g(\mu_{\eps_k}) dx \ \textrm{ for every Borel set } E \subset \Om.$$
As in Step 1, we may pass along a subsequence to a weak-$*$ limit 

\begin{equation}
\lambda << \Ln + |\mu_s|.
\label{UBAC}
\end{equation}
To prove the upper bound, we will show that

\begin{equation}
\label{GUBAC}
\frac{d \lambda}{d \Ln} (x_0) \leq g( \mu_{ac} (x_0)) \textrm{ for } \Ln a.e.  \, x_0 \in \Omega,
\end{equation}
and
\begin{equation}
\label{GUBS}
\frac{d \lambda}{d |\mu_s|} (x_0) \leq g^{\infty}( \nu(x_0)) \textrm{ for } |\mu_s| \, a.e. \, x_0 \in \Omega.
\end{equation}
Assuming ($\ref{GUBAC}$) and ($\ref{GUBS}$), by our boundary regularity assumption on $\Om$ and (\ref{UBAC}) we have $\lambda( \partial \Om) = 0$, and therefore

$$
\lambda(\Om) = \lambda(\overline{\Om}) \geq \limsup_{k \rga \infty} \lambda_k (\overline{\Om}) = \limsup_{k \rga \infty} \lambda_k(\Om),
$$
while

$$
\lambda(\Om) \leq \int_{\Om} g(\mu_{ac}(x)) dx + \int_{\Om} g^{\infty}(\nu(x)) d |\mu_s|(x),
$$
and putting these together, we have ($\ref{lemUB}$).
To prove ($\ref{GUBAC}$), we know that for $\Ln$-almost every $x_0$, we have

$$
\frac{ d \lambda}{d \Ln} (x_0) = \lim_{r \rga 0} \frac{\lambda(\overline{Q(x_0, r)})}{r^N},
$$

\begin{equation}
\lim_{r \rga 0} \frac{|\mu_s|(\overline{Q(x_0, r)})}{r^N}  = 0,
\label{singsmall}
\end{equation}
and
\begin{equation}
\lim_{r \rga 0}  \intav_{Q(x_0, r)} | \mu_{ac}(x) - \mu_{ac}(x_0) | dx = 0.
\label{acint}
\end{equation}
For all such points $x_0$, we can pick a sequence $r_n \rga 0$ such that $\lambda(\partial Q(x_0, r_n)) = 0$. Then,

$$
\frac{d \lambda}{d \Ln}(x_0) = \lim_{n \rga \infty} \frac{\lambda(\overline{Q(x_0, r_n)})}{r_n^N}  = \lim_{n \rga \infty} \lim_{k \rga \infty} \frac{\lambda_k(\overline{Q(x_0, r_n)}) }{r_n^N},
$$
and thus

\begin{align*}
\frac{d \lambda}{d \Ln} (x_0) &= \lim_{n \rga \infty} \lim_{k \rga \infty} \frac{1}{r_n^N} \int_{Q(x_0, r_n)} g( \mu * \phi_{\eps_k} ) dx\\
&= \lim_{n \rga \infty} \lim_{k \rga \infty} \frac{1}{r_n^N} \int_{Q(x_0, r_n)} g( \mu_{ac} * \phi_{\eps_k}  + \mu_{s} * \phi_{\eps_k} ) dx.
\end{align*}
By convexity of $g$, for any $p, q \in \Rd$ and $\theta \in (0,1)$ we have 

$$g(p + q) \leq \theta g \bigg(\frac{1}{\theta} p \bigg) + (1 - \theta) g \bigg( \frac{1}{1-\theta} q \bigg),$$
hence

$$
\frac{d \lambda}{d \Ln} (x_0) \leq  \varliminf_{n \rga \infty} \varliminf_{k \rga \infty} \frac{1}{r_n^N} \int_{Q(x_0, r_n)} \Bigg[ \theta g \bigg( \frac{1}{\theta} \mu_{ac} * \phi_{\eps_k} \bigg) + (1 - \theta) g \bigg( \frac{1}{1-\theta} \mu_s * \phi_{\eps_k} \bigg) \Bigg] dx
$$
and in view of Lemma $\ref{convolutions}$ we have

\begin{equation}
\frac{d \lambda}{d \Ln} (x_0) \leq \varliminf_{n \rga \infty} \varliminf_{k \rga \infty} \frac{1}{r_n^N} \int_{Q(x_0, r_n)} \Bigg[ \theta g\bigg(\frac{1}{\theta} \mu_{ac}\bigg) * \phi_{\eps_k} + (1 - \theta) g \bigg( \frac{1}{1-\theta} \mu_s * \phi_{\eps_k} \bigg) \Bigg] dx.
\label{i}
\end{equation}
Moreover,

\begin{align*}
\frac{1- \theta}{r_n^N} \int_{Q(x_0, r_n)}  g \bigg( \frac{1}{1-\theta} \mu_s * \phi_{\eps_k} \bigg) dx &\leq \frac{1- \theta}{r_n^N} \int_{Q(x_0, r_n)} C \bigg( 1 + \bigg|\frac{1}{1-\theta} \mu_s * \phi_{\eps_k}\bigg| \bigg) dx\\
&\leq (1 - \theta) C + \frac{C}{r_n^N} \int_{Q(x_0, r_n)} |\mu_s * \phi_{\eps_k}| dx,  
\end{align*}
and this yields

\begin{align}
\notag
\varlimsup_{n \rga \infty} \varlimsup_{k \rga \infty} \frac{1 - \theta}{r_n^N} \int_{Q(x_0,r_n)} g \bigg( \frac{1}{1-\theta} \mu_s * \phi_{\eps_k} \bigg) dx &\leq (1 - \theta) C + C \lim_{n \rga \infty} \frac{|\mu_s|(\overline{Q(x_0, r_n)})}{r_n^N}\\
&= (1-\theta) C,
\label{ii}
\end{align}
where we have used ($\ref{singsmall}$) and, from \cite{AFP00}, Theorem 2.2 that  
$$\limsup_{k \rga \infty} \int_{Q(x_0, r)} |\mu_s * \phi_{\eps_k}| dx \leq \limsup_{k \rga \infty} |\mu_s|(Q(x_0, r_n + \eps_k)) = |\mu_s|( \overline{Q(x_0, r_n)}).$$
In turn,

\begin{align}
\notag
\lim_{n \rga \infty} \lim_{k \rga \infty} \frac{1}{r_n^N} \int_{Q(x_0, r_n)} \theta g\bigg(\frac{1}{\theta} \mu_{ac}\bigg) * \phi_{\eps_k} dx &= \lim_{n \rga \infty} \frac{1}{r_n^N} \int_{Q(x_0, r_n)} \theta g\bigg( \frac{1}{\theta} \mu_{ac} \bigg) dx\\
\label{iii}
&= \theta g \bigg( \frac{1}{\theta} \mu_{ac} (x_0) \bigg).
\end{align}
To see why the last step above holds, note $g$ is convex with linear growth, and so by regularity properties of convex functions (see \cite{FL07}, Proposition 4.64) we have that $g$ is Lipschitz with some constant $L>0$. Thus,

$$
\bigg| \frac{1}{r^N} \int_{Q(x_0, r)} \bigg[ g \bigg( \frac{1}{\theta} \mu_{ac}(x) \bigg) - g \bigg( \frac{1}{\theta} \mu_{ac}(x_0) \bigg) \bigg] dx \bigg| \leq \frac{1}{r^N} \int_{Q(x_0, r)} \frac{L}{\theta} |\mu_{ac}(x) - \mu_{ac}(x_0)| dx \buildrel{r \rga 0}\over\longrightarrow 0
$$
\noindent
by virtue of ($\ref{acint}$).
By ($\ref{i}$), ($\ref{ii}$), and ($\ref{iii}$), we have for every $\theta \in (0,1)$ that

$$
\frac{d \lambda}{d \Ln} (x_0) \leq \theta g \bigg( \frac{1}{\theta} \mu_{ac}(x_0) \bigg) + (1-\theta) C,
$$
and letting $\theta \rga 1^-$, by continuity of $g$ we conclude that

$$
\frac{d \lambda}{d \Ln} (x_0) \leq g( \mu_{ac}(x_0))
$$
for $\Ln$ almost every $x_0 \in \Om$.

\

Next, we tackle the singular part, i.e., ($\ref{GUBS}$) . We know that for $|\mu_s|$ almost every $x_0$, we have

$$
\frac{ d \lambda}{d |\mu_s|} (x_0) = \lim_{r \rga 0} \frac{ \lambda( \overline{Q(x_0, r)})}{ |\mu_s|(\overline{Q(x_0, r)})},
$$

\begin{equation}
\lim_{r \rga 0} \frac{ r^N}{|\mu_s|(\overline{Q(x_0, r)})} = 0,
\label{lebsmall}
\end{equation}

\begin{equation}
\lim_{r \rga 0} \frac{1}{|\mu_s|(\overline{Q(x_0, r)})} \int_{Q(x_0, r)} |\mu_{ac}|(x) dx = 0,
\label{acsmall}
\end{equation}
and

$$
\lim_{r \rga 0} \intav_{\overline{Q(x_0, r)}} g^{\infty}(\nu(x)) d|\mu_s|(x) = g^{\infty}(\nu(x_0)).
$$
We choose a sequence $r_n \rga 0$ such that $|\mu_s|(\partial Q(x_0, r_n)) = 0.$ Note that by $(\ref{UBAC})$ we also have $\lambda(\partial Q(x_0, r_n)) = 0.$ We obtain

\begin{align*}
\frac{ d \lambda}{d |\mu_s|} (x_0) &= \lim_{n \rga \infty} \frac{ \lambda( Q(x_0, r_n))}{ |\mu_s|(Q(x_0, r_n))} = \lim_{n \rga \infty} \lim_{k \rga \infty} \frac{ \lambda_k(Q(x_0, r_n))}{ |\mu_s|(Q(x_0, r_n))}\\
&= \lim_{n \rga \infty} \lim_{k \rga \infty} \frac{1}{|\mu_s|(Q(x_0, r_n))} \int_{Q(x_0, r_n)} g( \mu * \phi_{\eps_k} ) dx.
\end{align*}
Again appealing to convexity of $g$, we get

\begin{equation}
\label{convexhalf}
|\mu_s| (x_0) \leq  \varliminf_{n \rga \infty} \varliminf_{k \rga \infty} \frac{1}{|\mu_s|(Q(x_0, r_n))} \int_{Q(x_0, r_n)} \bigg[ \frac{1}{2} g ( 2 \mu_{ac} * \phi_{\eps_k} ) + \frac{1}{2} g ( 2 \mu_s * \phi_{\eps_k} ) \bigg] dx.
\end{equation}
Since by (\ref{A1})

$$
 \frac{1}{2} g ( 2 \mu_{ac} * \phi_{\eps_k} ) \leq C( 1 + |\mu_{ac} * \phi_{\eps_k}| ) ,
$$
we have 

\begin{align}
\notag
\varlimsup_{n \rga \infty} \varlimsup_{k \rga \infty} &\frac{1}{|\mu_s|(Q(x_0, r_n))} \int_{Q(x_0, r_n)} \frac{1}{2} g ( 2 \mu_{ac} * \phi_{\eps_k} ) dx\\
\notag
\leq 
&\varlimsup_{n \rga \infty} \varlimsup_{k \rga \infty} \frac{1}{|\mu_s|(Q(x_0, r_n))} \int_{Q(x_0, r_n)} C( 1 + |\mu_{ac} * \phi_{\eps_k}| ) dx \\
\label{p} 
=&\lim_{n \rga \infty} \frac{C}{|\mu_s|(Q(x_0, r_n))}\bigg(  r_n^N + \int_{Q(x_0, r_n)} |\mu_{ac}|(x) dx \bigg) = 0,
\end{align}
where we used ($\ref{lebsmall}$) and ($\ref{acsmall}$).

Next we restrict our attention to the singular part in (\ref{convexhalf}). Noting that
\[ |(\phi_{\eps_k} * \mu_s)(x)| \leq (\phi_{\eps_k} * |\mu_s|)(x) =: t_{\eps_k} (x), \]
we have

\begin{align}
\int_{Q(x_0,r_n) \cap \{ t_{\eps_k} < 1 \}} g(2 \phi_{\eps_k} * \mu_s(x) ) dx \leq
r_n^N  \| g\|_{L^{\infty}(B(0,2))} = C r_n^N.
\label{o}
\end{align}
\noindent
Meanwhile, in the second region, using Lemma $\ref{convolutions}$ we obtain 

\begin{align}
\notag
\int_{Q(x_0,r_n) \cap \{ t_{\eps_k} \geq 1 \}} & g(2 \phi_{\eps_k} * \mu_s(x) ) dx\\
\leq &\int_{Q(x_0,r_n) \cap \{ t_{\eps_k} \geq 1 \}} \int_{\Om} \frac{g(2 t_{\eps_k}(x) \nu(y))}{t_{\eps_k}(x)}  \phi_{\eps_k}(x-y) d|\mu_s|(y)  dx.
\label{tbig}
\end{align}
Note that by (\ref{A3}),

$$
\bigg|\frac{g(2t_{\eps_k}(x) \nu(y))}{t_{\eps_k}(x)} - 2g^{\infty}(\nu(y)) \bigg| \leq \frac{C}{t_{\eps_k}(x)^{\alpha}},
$$
and by ($\ref{tbig}$) we have

\begin{align}
\notag
\int_{Q(x_0,r_n) \cap \{ t_{\eps_k} \geq 1 \}} g(2 \phi_{\eps_k} * \mu_s(x) ) &\leq \int_{Q(x_0,r_n) \cap \{ t_{\eps_k} \geq 1 \}} \int_{\Om} 2 g^{\infty}(\nu(y)) \phi_{\eps_k}(x-y) d|\mu_s|(y) dx\\
\notag
& \ \ \  + C\int_{Q(x_0,r_n) \cap \{ t_{\eps_k} \geq 1 \}} \int_{\Om} \frac{1}{t_{\eps_k}(x)^{\alpha}} \phi_{\eps_k}(x-y) d|\mu_s|(y) dx\\
\notag
&\leq \int_{Q(x_0,r_n)} 2 \big( g^{\infty}(\nu(\cdot))|\mu_s| * \phi_{\eps_k} \big) (x) dx  \\
\notag
& \ \ \  + C \int_{Q(x_0,r_n) \cap \{ t_{\eps_k} \geq 1 \}} \frac{1}{t_{\eps_k}(x)^{\alpha-1}} dx\\
&\leq \int_{Q(x_0,r_n)} 2\big( g^{\infty}(\nu(\cdot))|\mu_s| * \phi_{\eps_k} \big) (x) dx + C r^N,
\label{oo}
\end{align}

In view of ($\ref{o}$), ($\ref{oo}$), we have shown for every $\theta$ that

$$\int_{Q(x_0,r_n)} g(2 \phi_{\eps_k} * \mu_s(x) ) dx \leq  \int_{Q(x_0,r_n)} 2 \big( g^{\infty}(\nu(\cdot))|\mu_s| * \phi_{\eps_k} \big) (x) dx + C r_n^N,
$$
therefore, 

\begin{align}
\notag
\frac{1}{|\mu_s|(Q(x_0, r_n))} &\varlimsup_{k \rga \infty} \int_{Q(x_0,r_n)} \frac{1}{2} g(2 \phi_{\eps_k} * \mu_s(x) ) dx \\
&\leq \intav_{Q(x_0, r_n)} g^{\infty}( \nu(x)) d|\mu_s|(x) +C \frac{r_n^N}{|\mu_s|(Q(x_0, r_n))},
\label{pp}
\end{align}

and by ($\ref{lebsmall}$), ($\ref{p}$), ($\ref{pp}$) we conclude that
 
$$
\frac{d \lambda}{d |\mu_s|} (x_0) \leq g^{\infty}(\nu(x_0)).
$$

\qed

\noindent\makebox[\linewidth]{\rule{\paperwidth}{0.4pt}}

In our application, we are interested in the case when the measure $\mu$ is the Hessian of a $BH$ function. It should be noted that the first order case, when $\mu$ is the gradient of some $BV$ function, an area-strict density theorem follows from the integral representation results of Fonseca and M\"uller \cite{FM93}, and Ambrosio and Dal Maso \cite{ADM92}, with no regularity assumption on the boundary.

To apply Theorem \ref{density} to a given $u \in BH(\Om; \Rd)$, the main obstacle is finding an extension of $u$ to a larger set $U$ such that $|D(\nabla u)|(\partial \Om) = 0$. In order to achieve a fairly general class of domains, we shall borrow from the construction of Stein \cite{S70}.

First we will define the extension in the case where $\Om$ is of type special Lipschitz. Recall that we say a set $\Om \subset \mathbb{R}^{N}$ is \textit{special Lipschitz} if there is a Lipschitz function $f: \Rn \rga \mathbb{R}$ such that

\be{sLip}
\Om = \{(x',x_N) \in \mathbb{R}^{N}: x_N > f(x') \}
\ee
where we use the notation $(x',x_N)$ to identify $\mathbb{R}^N$ with $\mathbb{R}^{N-1} \times \mathbb{R}$.

We begin with a simpler approximation lemma.

\bl{translation}

Let $\Om \subset \mathbb{R}^{N}$ be a special Lipschitz domain. For any $u \in BH(\Om; \Rd)$ there exists a sequence $\{u_n\} \in W^{2,1}(\Om; \Rd)$ such that 

$$\lim_{n \rga \infty} \|u_n - u \|_{W^{1,1}(\Om; \mathbb{R}^d)} = 0 \textit{ and } \sup_{n} \|u_n\|_{W^{1,2}(\Om; \mathbb{R}^d)} < \infty.$$

\el

\proof

Consider a function $f$ as in (\ref{sLip}). Given any $v \in L^1(\Om)$ and $\delta > 0$, we define its translation $T_{\delta}v \in L^1(\Om_{\delta})$ via

$$T_{\delta}v(x) := v(x', x_N + \delta) $$
\noindent
where $\Om_{\delta} := \{(x',x_N) \in \mathbb{R}^{N}: x_N > f(x') - \delta \}$.

Note that if $L$ is the Lipschitz constant of $f$, then for $x \in \Om$ and $\eps < \frac{\delta}{1+L}$ we have $B(x, \eps) \subset \Om_{\delta}$, so the function $\phi_{\eps} * T_{\delta}v \in C^{\infty}(\Om)$ is well-defined.
Since $\nabla T_{\delta}u = T_{\delta} \nabla u$ and the translation is continuous in the $L^1$ norm, we have

$$\lim_{\delta \rga 0^+} \| T_{\delta}u - u \|_{W^{1,1}(\Om; \mathbb{R}^d)}  = 0.$$
By standard mollification results, we have for every $\delta>0$

$$\lim_{\eps \rga 0^+} \| \phi_{\eps} * T_{\delta} u - T_{\delta}u\|_{W^{1,1}(\Om; \mathbb{R}^d)} = 0, $$

$$ \limsup_{\eps \rga 0^+} \| \phi_{\eps} * T_{\delta} u \|_{W^{2,1}(\Om; \mathbb{R}^d)} \leq \|u \|_{BH(\Om)}.$$
This follows from the general fact that for any Radon measure $\mu \in \mathcal{M}(\Om_{\delta})$, if $\eps < \frac{\delta}{1+L}$ we have 

$$\int_{\Om_{\delta}} |\phi_{\eps} * \mu| dx \leq |\mu|(\Om).$$
Indeed,
\begin{align}
\notag
\int_{\Om_{\delta}} |\phi_{\eps} * \mu| dx &\leq \int_{\Om_{\delta}} \int_{\mathbb{R}^N} \phi_{\eps}(x-y) d|\mu|(y) dx  = \int_{\mathbb{R}^N} \int_{\Om_{\delta}} \phi_{\eps}(x-y) dx \ d|\mu|(y)  \\
&=  \int_{\mathbb{R}^N} \int_{\Om_{\delta} \cap B(y, \eps)} \phi_{\eps}(x-y) dx \ d|\mu|(y) =  \int_{\mathbb{R}^N} \int_{(\Om_{\delta} - y) \cap B(0, \eps)} \phi_{\eps}(z) dz d|\mu|(y)
\label{convdelt}
\end{align}
using the substitution $z = x-y$. Now, if the set $(\Om_{\delta} - y) \cap B(0, \eps)$ is nonempty, we must have dist$(y, \Om_{\delta}) < \eps$, and so $y \in  \Om$ because $\eps < \delta$. Thus we can bound (\ref{convdelt}) by
$$\int_{\Om} \int_{B(0,\eps)} \phi_{\eps}(z) dz d|\mu|(y) = |\mu|(\Om). $$
Applying this to $\phi_{\eps} * T_{\delta}u$, we see that
$$\int_{\Om} |\nabla^2(\phi_{\eps}*T_{\delta}u)| dx =  \int_{\Om} |\phi_{\eps} * D^2 (T_{\delta}u)| dx = \int_{\Om_{\delta}} |\phi_{\eps} * D^2 u| dx \leq |D^2u|(\Om)$$
and similarly for $\phi_{\eps} * T_{\delta} u$ and $\nabla(\phi_{\eps} * T_{\delta} u)$.

Thus, for any sequence $\delta_n \rga 0^+$ we can choose $\eps_n < \frac{\delta_n}{1+L}$ such that the smooth (and thus $W^{2,1}$) functions $\phi_{\eps_n} * T_{\delta_n} u$ converge to $u$ with bounded $W^{2,1}$ norm.

\qed

\bt{densitySL}
\label{densitySL}
Let $\Om \subset \mathbb{R}^{N}$ be a special Lipschitz domain. For any function $u \in BH(\Om; \Rd)$ there is an extension $E[u] \in BH(\mathbb{R}^{N}; \Rd)$ such that $E[u] = u$ in $\Om$ and $|D(\nabla E[u])|(\partial \Om) = 0.$

\et

The general theory of extending $BH$ functions can be reduced to the theory of extending $BV$ functions. We recall that (see \cite{EG92}, 5.4.1) since $\Om$ is Lipschitz, given $w_1 \in BV(\Om)$ and $w_2 \in BV(\mathbb{R}^{N} \setminus \overline{\Om})$, the function 

\[ w(x) := \begin{cases} 
      w_1(x) & \textrm{ for $x \in \Om$}, \\
      w_2(x) & \textrm{ for $x \in \mathbb{R}^{N} \setminus \overline{\Om}$}, \\
   \end{cases} \]
is a $BV$ function with

\begin{equation}
Dw = Dw_1 \mres \Om + Dw_2 \mres (\mathbb{R}^{N} \setminus \overline{\Om})  + (\textrm{Trace}(w_2) - \textrm{Trace}(w_1)) \nu \ \mathcal{H}^{N-1} \mres \partial \Om
\label{BV trace}
\end{equation}
where $\nu$ is the outward normal vector to $\partial \Om$.

\

Let $u \in BH(\Om)$. Since $\nabla u$ is a $BV$ function, in view of (\ref{BV trace}) to guarentee that our extension does not charge the boundary, it suffices to ensure that the traces of $\nabla u$ and of the extension $\nabla E[u]$ agree on the boundary $\partial \Om$.

\

We will be inspired by the construction given by Stein to introduce $E[u]$, namely

\be{ext}
E[u](x) := \begin{cases}
\displaystyle \int_1^{\infty} u(x', x_N + \lambda \delta(x)) \psi(\lambda) d \lambda & \textrm{ for $x \in \mathbb{R}^{N} \setminus \overline{\Om}$}, \\
\\
\displaystyle u(x) & \textrm{ for $x \in \Om$}, \\
\end{cases}
\ee
where $\delta \in C^{\infty}(\mathbb{R}^{N})$ is a regularized distance function, that is, 

\begin{equation}
c \ \textrm{dist}(x, \partial \Om) \leq \delta(x) \leq C \ \textrm{dist}(x, \partial \Om),
\label{likedist}
\end{equation}
\begin{equation}
\bigg| \frac{\partial^{\alpha}}{\partial x^{\alpha}} \delta (x) \bigg| < C_{\alpha} \bigg(\frac{1}{\textrm{dist}(x, \partial \Om)}\bigg)^{1-|\alpha|},
\label{derivstoo}
\end{equation}
for every multi-index $\alpha$ and $\psi: [1, \infty) \rga \mathbb{R}$ is a continuous function with

\be{moments}
\int_1^{\infty} \psi(\lambda) d \lambda = 1, \ \int_1^{\infty} \psi(\lambda) \lambda^k d \lambda = 0, \ k \in \mathbb{N},
\ee
and
\be{decay}
\lambda^k \psi( \lambda) \rga 0 \textrm{ as } \lambda \rga + \infty \textrm{ for all } k \in \mathbb{N}.
\ee

In fact we will use a different regularized distance from the one used by Stein. In general the distance function is a Lipschitz function, and thus will be differentiable almost everywhere. However, the gradient of the regularized distance of Stein does not in general approach the gradient of the distance function as a pointwise limit. We elect to use a regularized distance introduced by Lieberman \cite{L85}, defined as follows.

\bd{Regularized Distance}
Let $\Om$ be an open subset of $\Rn$. We define the signed distance to $\partial \Om$ by

\[ d(x) := \begin{cases} 
      \textnormal{dist}(x, \partial \Om) & \textrm{ for $x \not \in \Om$}, \\
      -\textnormal{dist}(x, \partial \Om)& \textrm{ for $x \in \Om$}. \\
   \end{cases} \]
For $x \in \Rn, t \in \mathbb{R}$, define
$$G(x,  t) := \int_{|z|<1} d\bigg(x-\frac{t}{2}z\bigg) \phi(z) dz$$
where $\phi = \phi_1$ is a mollifier with $\textnormal{supp}(\phi) \subset \subset B(0,1)$.
Then a regularized distance is given by the equation
\begin{equation}
\label{rhodef}
\rho(x) := G(x, \rho(x)).
\end{equation}
\ed

The regularized distance of Lieberman satisfies the necessary bounds (\ref{likedist}), (\ref{derivstoo}) as noted in \cite{L85} Lemma 1.1, and hence it may be used as a substitute for the regularized distance of Stein. In fact, in the construction of Stein, the function $\delta$ is not the regularized distance, but a scaled multiple of it, where this scalar factor, depending on the Lipschitz constant of the domain, is chosen so that for every $x \not \in \Om$ we have (see (\ref{ext}))
$$x + \delta(x) e_N \in \Om.$$
In light of this freedom, we will make this factor slightly larger in order to guarantee the bound
\begin{equation}
\label{Nbound}
\frac{\partial \delta}{\partial x_N}(x) \leq -2,
\end{equation}
and this is possible because it can be shown (see Lemma \ref{Lieb}) that there exists a constant $C = C(\Om) > 0$ such that
$$\frac{\partial \rho}{\partial x_N}(x) \leq -C$$ for almost every $x \in \mathbb{R}^N$.

\bl{Lieb}
Let $\Om \subset \Rn$ be a special Lipschitz domain with constant $L$. Then for every $x \in \mathbb{R}^N \setminus \overline{\Om}$
$$\frac{\partial \rho}{\partial x_N}(x) \leq \frac{-2}{3} \frac{1}{1+L^2}. $$
\el

First we note that since $\Om$ is a special Lipschitz domain, for every $x \in \mathbb{R}^N \setminus \overline{\Om}$ and $h >0$ sufficiently small so that $x + h e_N \in \mathbb{R}^N \setminus \overline{\Om}$, we have
\be{out}
d(x + h e_N) \leq d(x) -  \frac{h}{1+L^2},
\ee
where $L$ is the Lipschitz constant of $f$. Indeed, let $x \in \mathbb{R}^N \setminus \overline{\Om}$ and $h>0$ be as above. Let  $y \in \partial \Om$ be arbitrary. We claim that
\be{anyY}
d(x+h e_N) \leq |x - y| - \frac{h}{1+L^2}.
\ee
If $y' = x'$, then $x, x + h e_N$ and $y$ are colinear and 
$$d(x + h e_N) \leq |x + h e_n - y| = |x-y| - h \leq  |x - y| - \frac{h}{1+L^2}. $$
Thus it suffices to prove (\ref{anyY}) in the case where $x$, $x + h e_N$, and $y$ are not colinear. Since $y \in \partial \Om$ and $\Om$ is special Lipschitz, we have the exterior cone condition
\be{cone}
K = \{ z \in \mathbb{R}^N : z_N \geq f(y') + L |z' - y'| \} \subset \overline{\Om}
\ee
as $z \not \in \overline{\Om}$ implies $z_N < f(z') \leq f(y') + L|z'-y'|$. Noting that the point $y+h e_N \in \Om$ is inside $K$ and the point $x + h e_N$ is outside of $K$, there is some $z \in \partial K$ on the line from $x+h e_N$ to $y + h e_N$. Since $z \in \overline{\Om}$, there is a point $\tilde{z} \in \partial \Om$ between $x + h e_N$ and $z$ such that $z, \tilde{z}, x+he_N$ are colinear. Thus,
\be{nuffZ}
d(x+ h e_N) \leq |x+ he_N - \tilde{z}| \leq |x + h e_N - z|.
\ee
As $x + he_N, z, y+h_N$ are colinear we have
\be{colin}
|x + h e_N - z| = |x+h e_N - y - h e_N| - |y + h e_N - z| = |x-y| - |y + h e_N - z|.
\ee
We have $z_N = y_N + L|z'-y'|$ because $z \in \partial K$ and $y \in \partial \Om$, and so 
\be{quad}
|y + h e_N - z|^2 = |z'-y'|^2 + (h - L |z'-y'|)^2.
\ee
Consider the quadratic function
$$ t \mapsto t^2 + (h-Lt)^2 $$
and observe that it is minimized at $t =t^* := \frac{hL}{1+L^2}.$ By (\ref{quad}) we have 
$$|y + h e_N - z|^2 \geq (t^*)^2 + (h-Lt^*)^2 \geq (h-Lt^*)^2= h^2 \bigg( \frac{1}{1+L^2} \bigg)^2,$$
i.e.,
$$|y + h e_N - z| \geq \frac{h}{1+L^2}. $$
Taking this into account, (\ref{colin}) yields
$$|x + h e_N - z| \leq |x-y| - \frac{h}{1+L^2}$$
which, with (\ref{nuffZ}), implies
$$d(x+ h e_N) \leq |x-y| - \frac{h}{1+L^2}.$$
Taking the infimum over all $y \in \partial \Om$, we deduce that
$$d(x + h e_N) \leq d(x) - \frac{h}{1+L^2}.$$
By a reflection argument, for every $x \in \Om, h >0$ sufficiently small so that $x - he_N \in \Om$, we have
\be{in}
d(x - he_N) \geq d(x) + \frac{h}{1+L^2}.
\ee
We consider the special Lipschitz set $U$ whose boundary is given by the function $g(x') := -f(-x')$. Note that $x_N < g(x')$ if and only if $-x_N > f(-x')$. Thus, $U = - (\mathbb{R}^N \setminus \overline{\Om})$ and $\mathbb{R}^N \setminus \overline{U} = - \Om$, and so the signed distance function $\tilde{d}$ for $U$ satisfies $\tilde d(x) = -d(-x)$. Thus for all $x \in \Om$, $h>0$ sufficiently small, by (\ref{in}) applied to $\tilde{d}$ at $-x$ we have
$$\tilde{d}(-x + h e_N)) \leq \tilde{d}(-x) - \frac{h}{1+L^2},$$ 
and thereby
$$-d(x -h e_N) \leq -d(x) - \frac{h}{1+L^2},  $$
or, equivalently,
$$d(x-h e_N) \geq d(x) + \frac{h}{1+L^2}.$$
By (\ref{in}), (\ref{out}), and the fact that $|\partial \Om| = 0$, we conclude that
\be{undLip}
\frac{\partial d}{\partial x_N}(x) \leq \frac{-1}{1+L^2}
\ee
for almost every $x \in \Rn$. As $d \in W^{1, \infty}$, $\phi \in C^{\infty}_c$, we can differentiate under the integral sign to get
$$ \frac{\partial G}{\partial t} (x,t) = \int_{|z|<1} - \nabla d \bigg(x- \frac{t}{2} z \bigg) \cdot \frac{z}{2} \phi(z) dz$$
and thus, since $\|\nabla d\|_{\infty} \leq 1$ and $\int_{B(0,1)} \phi = 1$, we obtain
\begin{equation}
\label{tder}
\bigg| \frac{\partial G}{\partial t} \bigg| \leq \frac{1}{2}.
\end{equation}
On the other hand, we have
$$\frac{\partial G}{\partial x_N} (x,t) = \int_{|z|<1} \frac{\partial d}{\partial x_N} \bigg(x-\frac{t}{2}z \bigg) \phi(z) dz$$
and so by (\ref{undLip}), 
\begin{equation}
\label{Nder}
\frac{\partial G}{\partial x_N} (x,t) \leq \frac{-1}{1+L^2} .
\end{equation}
We differentiate (\ref{rhodef}) to obtain
$$\frac{\partial \rho}{\partial x_N} (x) =  \frac{\partial G}{\partial x_N}(x, \rho(x)) + \frac{\partial G}{\partial t}(x, \rho(x)) \frac{\partial \rho}{\partial x_N} (x) . $$
Solving this equation for $\frac{\partial \rho}{\partial x_N}$ yields
$$\frac{\partial \rho}{\partial x_N} (x) = \frac{ \frac{\partial G}{\partial x_N}(x, \rho(x))}{1-\frac{\partial G}{\partial t}(x, \rho(x))} \leq - \bigg( \frac{2}{3} \bigg)  \frac{1}{1+L^2}$$
where we have used the bounds (\ref{tder}) and (\ref{Nder}). 
\qed

We now have the tools to prove Theorem \ref{densitySL}.  \\
\noindent
\textit{Proof of Theorem \ref{densitySL}.} We claim that the function
\be{}
E[u](x) := \begin{cases}
\displaystyle \int_1^{\infty} u(x', x_N + \lambda \delta(x)) \psi(\lambda) d \lambda & \textrm{ for $x \in \mathbb{R}^{N} \setminus \overline{\Om}$}, \\
\\
\displaystyle u(x) & \textrm{ for $x \in \Om$}, \\
\end{cases}
\ee
is our desired extension. First, we prove that $E[u] \in BH(\mathbb{R}^{N}; \Rd)$. Consider an approximating sequence $\{u_n \} \subset W^{2,1}(\Om; \Rd)$ as in Lemma \ref{translation}. Since $\{ u_n \}$ is bounded in $W^{2,1}(\Om; \Rd)$ and the extension operator $E$ is a continuous linear operator

$$E : W^{2,1}(\Om; \Rd) \rga W^{2,1}(\mathbb{R}^{N}; \Rd), $$
the sequence $\{ E[u_n] \}$ is bounded in $W^{2,1}(\mathbb{R}^{N}; \Rd)$, and thus, along a subsequence (not relabeled), there is a function $v \in BH(\mathbb{R}^{N}; \Rd)$ such that

$$ E[u_n] \rga v \textrm{ in } W^{1,1}(\Rn; \Rd).$$
We claim that $v = E[u]$. To see this, first take $(x',x_N) \in \mathbb{R}^{N} \setminus \overline{\Om}$. Then,

\begin{align*}
\big|E[u](x) - E[u_n](x)\big| &\leq \int_1^{\infty} C |u(x', x_N + \lambda \delta(x)) - u_n(x',x_N + \lambda \delta(x))| d \lambda\\
&=\int_{\delta(x)}^{\infty} \frac{C}{\delta(x)} |u(x', x_N + s) - u_n(x',x_N + s)| d s.  
\end{align*}
For fixed $x_N$ and $\ell >0$, we can integrate both sides with respect to $x'$ over the set $\mathbb{R}^{N-1} \cap \{\delta(\cdot,x_N) \geq \ell \}$ to get

\begin{align*} &\int_{\mathbb{R}^{N-1} \cap \{\delta(\cdot,x_N) \geq \ell \}} |E[u](x',x_N) - E[u_n](x',x_N)| dx' \\
&\leq \frac{C}{\ell} \int_{\mathbb{R}^{N-1} \cap \{\delta(\cdot,x_N) \geq \ell \}} \int_{\ell}^{\infty} |u(x', x_N + s) - u_n(x',x_N + s)| d s dx'\\
&= \frac{C}{\ell} \int_{\mathbb{R}^{N-1} \cap \{\delta(\cdot,x_N) \geq \ell \}} \int_{\ell+x_N}^{\infty} |u(x', \tau) - u_n(x', \tau)| d \tau dx' \\
& \leq \frac{C}{\ell} \int_{\Om} |u(x',\tau) - u_n(x', \tau)| dx' d \tau.
\end{align*}
For any $T < \infty$ we integrate over $x_N \in (0, T)$ to get

$$\int_0^T \int_{\mathbb{R}^{N-1} \cap \{\delta(\cdot,x_N) \geq \ell \}} |E[u](x) - E[u_n](x)| dx \leq C \frac{T}{\ell} \int_{\Om} |u(x',\tau) - u_n(x', \tau)| dx' d \tau. $$ 

Since for every fixed $T$ and $\ell$, the right hand side goes to $0$ as $n \rga \infty$, we see that $\{ E[u_n] \}$ converges to $E[u]$ in $L^1_{loc}(\Rn \setminus \overline{\Om}; \Rd)$. However, we also know that $\{ E[u_n] \}$ converges in $W^{1,1}(\Rn; \Rd)$ to $v$, so we must have $E[u] = v$, and therefore $E[u] \in BH(\Rn; \Rd)$. We claim that 

$$\textrm{Trace}(E[u] ; \partial (\overline{\Om}^\mathsf{c})) = \textrm{Trace}(u; \partial \Om)$$
and
$$\textrm{Trace}(\nabla E[u] ; \partial (\overline{\Om}^\mathsf{c})) = \textrm{Trace}(\nabla u; \partial \Om).$$

To see this, we will use the following characterization of the trace operator, found in \cite{EG92} Theorem 5.3.2. Since $E[u]$ is a $BV$ function, for $\mathcal{H}^{N-1}$ almost every $x_0 \in \partial (\overline{\Om}^\mathsf{c})$, we have
\begin{align}
\notag
\textrm{Trace}(E[u] ; \partial (\overline{\Om}^\mathsf{c}))(x_0) &= \lim_{r \rga 0^+} \intav_{\Om^\mathsf{c} \cap B(x_0,r)} E[u](x) dx.\\
&= \lim_{r \rga 0^+} \intav_{\Om^\mathsf{c} \cap B(x_0,r)} \int_1^\infty u(F_{\lambda}(x)) \psi(\lambda) d \lambda d x  \\
&= \lim_{r \rga 0^+} \int_1^\infty \intav_{\Om^\mathsf{c} \cap B(x_0,r)} u(F_{\lambda}(x)) d x  \ \psi(\lambda)  d \lambda ,
\label{psiball}
\end{align}
where
$$F_{\lambda}(x) := x + \lambda \delta(x) e_N, \ x \in \overline{\Om}^\mathsf{c},$$
with
$$\nabla F_{\lambda}(x) = I + \lambda e_N \otimes \nabla \delta(x) $$
and
\be{dets}
\det(\nabla F_{\lambda}(x)) = 1 + \lambda \frac{\partial \delta}{\partial x_N}(x) \leq -1
\ee
where we have used (\ref{Nbound}) and the fact that $\lambda \geq 1$. In particular this implies that $F_{\lambda}$ is a local diffeomorphism. Thus we have $u \circ F_{\lambda} \in BV(\overline{\Om}^\mathsf{c}; \Rd)$ for every $\lambda \geq 1$, and
\begin{align}
\notag
\lim_{r \rga 0} \intav_{\overline{\Om}^\mathsf{c} \cap B(x_0,r)} u(F_{\lambda}(x)) d x &= \textrm{Trace}( u \circ F_{\lambda} ; \partial (\overline{\Om}^\mathsf{c}))(x_0) = \textrm{Trace}(u ; \partial \Om)(F_{\lambda}(x_0)) \\ &= \textrm{Trace}(u ; \partial \Om)(x_0)
\label{pwtrace}
\end{align}
for $\mathcal{H}^{N-1}$ almost every $x_0$ on $\partial \Om$, since $F_{\lambda}(x) = x$ on $\partial \Om$. We also observe that since $\Om^\mathsf{c}$ is locally a set of finite perimeter, for $\mathcal{H}^{N-1}$ almost every $x_0 \in \partial \Om$ we have there is $r_0$ such that for $r \in (0, r_0)$ we have
\be{finiteper}
\frac{\Om^\mathsf{c} \cap B(x_0, r)}{r^N} \geq \frac{1}{3}.
\ee

Now, for all $\lambda \geq 1$ we obtain 

\begin{align}
\notag
\bigg| \intav_{\Om^\mathsf{c} \cap B(x_0,r)} u(F_{\lambda}(x)) d x \bigg|&\leq \frac{1}{|\Om^\mathsf{c} \cap B(x_0,r)|} \int_{F_{\lambda}(\Om^\mathsf{c} \cap B(x_0,r))} \frac{|u(y)|}{|\det(\nabla F_{\lambda}(F_{\lambda}^{-1}(y)))|} dy 
\\
\notag
& \leq \frac{1}{|\Om^\mathsf{c} \cap B(x_0,r)|} \int_{\Om \cap B(F_{\lambda}(x_0), C \lambda r)} |u(y)| dy
\\
& = \frac{1}{|\Om^\mathsf{c} \cap B(x_0,r)|} \int_{\Om \cap B(x_0, C \lambda r)} |u(y)| dy
\label{L1bound}
\end{align}
for some constant $C(\Om)$, where we use (\ref{dets}) and the facts that since $\lambda \geq 1$
$$|F_{\lambda}(x)-F_{\lambda}(x_0)| = |x-x_0 + \lambda (\delta(x)-\delta(x_0)) e_N| \leq \lambda |x-x_0| + \lambda \textrm{Lip}(\delta) |x-x_0|, $$
which, since $F_{\lambda}(x_0) = x_0$, yields
$$F_{\lambda}(B(x_0, r)) \subset B(x_0, C \lambda r)$$
for some constant $C > 0$. Applying (\ref{finiteper}), we can estimate
\begin{equation}
\label{bigball}
\frac{|\Om \cap B(x_0,C \lambda r)|}{|\Om^\mathsf{c} \cap B(x_0,r)|} \leq \frac{|B(x_0, C \lambda r)|}{\frac{1}{3} r^N } = \frac{3 C^N \lambda^N r^N}{r^N} = C(\Om) \lambda^N.
\end{equation}
Now,
\[ \lim_{s \rga 0^+} \frac{1}{|\Om \cap B(x_0,s)|} \int_{\Om \cap B(x_0, s)} |u(y)| dy = \textrm{Trace}(|u|; \partial \Om)(x_0) \]
implies that there is some $s_0>0$ such that for $0<s<s_0$ we have
\[ \frac{1}{|\Om \cap B(x_0,s)|} \int_{\Om \cap B(x_0, s)} |u(y)| dy \leq 2 \textrm{Trace}(|u|; \partial \Om)(x_0). \]
Thus, for $\lambda < \frac{s_0}{Cr}$ we have
\[  \frac{1}{|\Om \cap B(x_0,C \lambda r)|} \int_{\Om \cap B(x_0, C \lambda r)} |u(y)| dy = 2 \textrm{Trace}(|u|; \partial \Om)(x_0) \]
and for $\lambda \geq \frac{s_0}{Cr} $ we have
\[ \frac{1}{|\Om \cap B(x_0,C \lambda r)|}  \int_{\Om \cap B(x_0, C \lambda r)} |u(y)| dy \leq  \frac{1}{|\Om \cap B(x_0,s_0)|} \int_{\Om} |u(y)| dy,  \]
and thus for every $\lambda \geq 1$ we have the bound
\begin{align}
\notag
\frac{1}{|\Om \cap B(x_0,C \lambda r)|}  \int_{\Om \cap B(x_0, C \lambda r)} |u(y)| dy &\leq \max \bigg\{2 \textrm{Trace}(|u|; \partial \Om)(x_0),   \frac{1}{|\Om \cap B(x_0,s_0)|} \int_{\Om} |u(y)| dy \bigg\} \\
&=: M.
\label{max}
\end{align}
By $(\ref{L1bound}), (\ref{bigball})$ and (\ref{max}), for every $\lambda \geq 1$, we have
\be{lambound}
\sup_{r \in (0, r_0) } \bigg| \intav_{\Om^\mathsf{c} \cap B(x_0,r)} u(F_{\lambda}(x)) d x \bigg| \leq C(\Om) \lambda^N M,
\ee
and so
\be{l1bound1}
\bigg| \intav_{\Om^\mathsf{c} \cap B(x_0,r)} u(F_{\lambda}(x)) d x \cdot \psi(\lambda) \bigg| \leq C(\Om) M \lambda^N |\psi(\lambda)| \in L^1([1, \infty))
\ee
for all $r \in (0, r_0)$, where we used (\ref{decay}). Thus we can apply the Dominated Convergence theorem in (\ref{psiball}) to conclude, by (\ref{pwtrace}), and we conclude that
\begin{align}
\notag
\textrm{Trace}(E[u]; \partial \overline{\Om}^\mathsf{c})(x_0)&= \lim_{r \rga 0^+} \int_1^\infty \intav_{\Om^\mathsf{c} \cap B(x_0,r)} u(F_{\lambda}(x)) d x \  \psi(\lambda)  d \lambda \\
&= \int_1^{\infty} \textrm{Trace}(u; \partial \Om)(x_0) \psi(\lambda) d \lambda = \textrm{Trace}(u; \partial \Om)(x_0),
\label{first}
\end{align}
where we have used (\ref{moments}). To see that the trace of $\nabla E[u]$ on $\partial \overline{\Om}^\mathsf{c}$ agrees with the trace of $\nabla u$ on $\partial \Om$, for $i \in \{ 1, ..., N \}$, note that

\be{dE}
\frac{\partial}{\partial x_i} E[u](x) = \int_{1}^{\infty} \frac{\partial u}{\partial x_i}(F_\lambda(x)) \psi(\lambda) d \lambda + \int_{1}^{\infty} \frac{\partial u}{\partial x_N}(F_\lambda(x)) \frac{\partial \delta}{\partial x_i}(x) \lambda \psi(\lambda) d\lambda.
\ee
Fix $x_0 \in \partial \Om$ such that Trace($\nabla u; \Om)(x_0)$ exists. Just as in (\ref{first}),
$$\lim_{r \rga 0^+} \intav_{\Om^\mathsf{c} \cap B(x_0,r)}  \int_{1}^{\infty} \frac{\partial u}{\partial x_i}(F_\lambda(x)) \psi(\lambda) d \lambda \ dx = \textrm{Trace}\bigg(\frac{\partial u}{\partial x_i}; \partial \Om\bigg)(x_0)$$
We claim that
\be{1o0}
\lim_{r \rga 0^+} \intav_{\Om^\mathsf{c} \cap B(x_0,r)} \int_{1}^{\infty} \frac{\partial u}{\partial x_N}(F_\lambda(x)) \frac{\partial \delta}{\partial x_i}(x) \lambda \psi(\lambda) d\lambda = 0,
\ee
Note that if (\ref{1o0}) holds, we conclude that
\begin{align*}
\textrm{Trace}\bigg( \frac{\partial E[u]}{\partial x_i} ; \partial \overline{\Om}^{\mathsf{c}} \bigg) &= \lim_{r \rga 0^+} \bigg( \intav_{\Om^\mathsf{c} \cap B(x_0, r)} \int_{1}^{\infty} \frac{\partial u}{\partial x_i}(F_\lambda(x)) \psi(\lambda) d \lambda + \int_{1}^{\infty} \frac{\partial u}{\partial x_N}(F_\lambda(x)) \frac{\partial \delta}{\partial x_i}(x) \lambda \psi(\lambda) d\lambda \bigg) \\
&= \textrm{Trace}\bigg(\frac{\partial u}{\partial x_i}; \partial \Om\bigg)(x_0).
\end{align*}
We will prove (\ref{1o0}) by showing that for every sequence $r_n \rga 0^+$, there is a subsequence $r_{n_k} \rga 0^+$ such that
$$\lim_{k \rga \infty} \intav_{\Om^\mathsf{c} \cap B(x_0,r_{n_k})} \int_{1}^{\infty} \frac{\partial u}{\partial x_N}(F_\lambda(x)) \frac{\partial \delta}{\partial x_i}(x) \lambda \psi(\lambda) d\lambda = 0. $$
By Fubini's theorem,
$$ \intav_{\Om^\mathsf{c} \cap B(x_0,r)} \int_{1}^{\infty} \frac{\partial u}{\partial x_N}(F_\lambda(x)) \frac{\partial \delta}{\partial x_i}(x) \lambda \psi(\lambda) d\lambda \ dx $$
$$= \int_{1}^{\infty} \intav_{\Om^\mathsf{c} \cap B(x_0,r)}  \frac{\partial u}{\partial x_N}(F_\lambda(x)) \frac{\partial \delta}{\partial x_i}(x) dx \ \lambda \psi(\lambda) d \lambda. $$
Now, for any sequence $r_n \rga \infty$, since $\frac{\partial \delta}{\partial x_i} \in L^{\infty}(\Rn)$ the sequence
$$\intav_{\Om^\mathsf{c} \cap B(x_0,r_n)}  \frac{\partial \delta}{\partial x_i}(x)   dx$$
is bounded and we can find a convergent subsequence such that
\be{bigD}
\lim_{k \rga \infty} \intav_{\Om^\mathsf{c} \cap B(x_0,r_{n_k})}   \frac{\partial \delta}{\partial x_i}(x) dx \textrm{ exists.}
\ee
We denote this limit $D(x_0)$. Now, again by the characterization of traces in \cite{EG92} Theorem 5.3.2, for all $\lambda \geq 1$ we have
$$\lim_{k \rga \infty} \intav_{\Om^\mathsf{c} \cap B(x_0, r_{n_k})} \bigg| \frac{\partial u}{\partial x_N}(F_\lambda(x)) - \textrm{Trace} \bigg( \frac{\partial u}{\partial x_N} ; \partial \Om \bigg)(x_0) \bigg| dx = 0.$$
Thus, for every $\lambda \geq 1$ we obtain
\begin{align*}
& \bigg| \intav_{\Om^\mathsf{c} \cap B(x_0, r_{n_k})} \frac{\partial u}{\partial x_N}(F_\lambda(x))  \frac{\partial \delta}{\partial x_i}(x) dx - D(x_0) \textrm{Trace} \bigg( \frac{\partial u}{\partial x_N} ; \partial \Om \bigg)(x_0) \bigg| \\
&\leq \bigg| \intav_{\Om^\mathsf{c} \cap B(x_0, r_{n_k})} \frac{\partial u}{\partial x_N}(F_\lambda(x))  \frac{\partial \delta}{\partial x_i}(x)  - \textrm{Trace} \bigg( \frac{\partial u}{\partial x_N} ; \partial \Om \bigg)(x_0)  \frac{\partial \delta}{\partial x_i}(x)  dx \bigg| \\
& \ \ \ + \bigg| \intav_{\Om^\mathsf{c} \cap B(x_0, r_{n_k})} \textrm{Trace} \bigg( \frac{\partial u}{\partial x_N} ; \partial \Om \bigg)(x_0)  \frac{\partial \delta}{\partial x_i}(x) - D(x_0) \textrm{Trace} \bigg( \frac{\partial u}{\partial x_N} ; \partial \Om \bigg)(x_0)  dx \bigg| \\
& \leq \bigg\| \frac{\partial \delta}{\partial x_i} \bigg\|_{\infty} \intav_{\Om^\mathsf{c} \cap B(x_0, r_{n_k})} \bigg| \frac{\partial u}{\partial x_N}(F_\lambda(x))   - \textrm{Trace} \bigg( \frac{\partial u}{\partial x_N} ; \partial \Om \bigg)(x_0) \bigg| dx \\
& \ \ \ + \bigg| \textrm{Trace} \bigg( \frac{\partial u}{\partial x_N} ; \partial \Om \bigg)(x_0) \bigg| \bigg| \intav_{\Om^\mathsf{c} \cap B(x_0, r_{n_k})}  \frac{\partial \delta}{\partial x_i}(x) dx - D(x_0) \bigg| .
\end{align*}
Letting $k \rga \infty$, using (\ref{bigD}) we conclude that 
$$\lim_{k \rga \infty} \intav_{\Om^\mathsf{c} \cap B(x_0,r_{n_k})}  \frac{\partial u}{\partial x_N}(F_\lambda(x)) \frac{\partial \delta}{\partial x_i}(x) dx = \textrm{Trace} \bigg( \frac{\partial u}{\partial x_N} ; \partial \Om \bigg)(x_0) D(x_0). $$
By the same argument as for (\ref{lambound}) we have
$$\sup_k \bigg|  \intav_{\Om^\mathsf{c} \cap B(x_0,r_{n_k})}  \frac{\partial u}{\partial x_N}(F_\lambda(x)) dx \bigg| \leq C(\Om)  \lambda^N  $$
and so it follows from $\frac{\partial \delta}{\partial x_i} \in L^{\infty}(\Rn)$ that
\be{l1bound2}
\lambda \psi(\lambda) \intav_{\Om^\mathsf{c} \cap B(x_0,r_{n_k})}  \frac{\partial u}{\partial x_N}(F_\lambda(x)) \frac{\partial \delta}{\partial x_i}(x) dx \in L^1([1,\infty)).
\ee
By the Dominated Convergence theorem, we have
\begin{align}
\notag
\lim_{k \rga \infty}   \int_{1}^{\infty} \intav_{\Om^\mathsf{c} \cap B(x_0,r_{n_k})}  \frac{\partial u}{\partial x_N}(F_\lambda(x)) \frac{\partial \delta}{\partial x_i}(x) dx \ \lambda \psi(\lambda) d \lambda \\
\label{moment2}
= \int_1^{\infty} \textrm{Trace}\bigg(\frac{\partial u}{\partial x_N} ; \partial \Om\bigg)(x_0) D(x_0) \lambda \psi(\lambda) d \lambda = 0,
\end{align}
where we use the fact that $\int_1^\infty \lambda \psi(\lambda) d \lambda = 0$. Thus in fact
$$\lim_{r \rga 0^+} \intav_{\Om^\mathsf{c} \cap B(x_0,r)} \int_{1}^{\infty} \frac{\partial u}{\partial x_N}(F_\lambda(x)) \frac{\partial \delta}{\partial x_i}(x) \lambda \psi(\lambda) d\lambda = 0. $$
We conclude that the traces of $E[u]$ and $\nabla E[u]$ on $\partial (\overline{\Om}^\mathsf{c})$ align with those of $u$ and $\nabla u$ on $\partial \Om$. Thus, we have our desired extension.

\qed

\br{} We remark that in the previous proof we did not use the full strength of properties (\ref{moments}) and (\ref{decay}). Indeed, (\ref{moments})$_1$ was needed in (\ref{first}) and (\ref{moments})$_2$ was used with $k=1$ in (\ref{moment2}). In (\ref{l1bound1}) it sufficed to have $\lambda^{N+2}\psi$ bounded, while  (\ref{l1bound2}) would hold with $\lambda^{N+3}\psi$ bounded.
\er

\bt{densityL}
\label{densityL}
Let $\Om \subset \Rn$ be an open, bounded Lipschitz set. For any function $u \in BH(\Om; \Rd)$ there exists an extension $E[u] \in BH(\Rn; \Rd)$ such that $|D(\nabla E[u])|(\partial \Om) = 0$.
\et

\proof

Since $\Om$ is a Lipschitz domain, we can cover $\Om$ by bounded open $U_0 \subset \subset \Om$ and $U_1, \dots , U_k$ such that $U_i \cap \partial \Om$ is the graph of a Lipschitz function. We may also choose a smooth partition of unity $\psi_0, \dots , \psi_k$ suboordinate to this cover. 

\

For $i \geq 1$, the domains $U_i$ are the subgraphs of Lipschitz functions, so we can find special Lipschitz domains $\Om_i$ such that $\Om_i \cap U_i = U_i \cap \Om$. Thus, by extending the functions $\psi_i u$ from $U_i \cap \Om$ to $\Om_i$ by zero, we can consider them to be defined on the special Lipschitz domains $\Om_i$. By Lemma \ref{densitySL}, we can find $BH$ functions $E[ \psi_i u] \in BH(\Rn; \Rd) $ which satisfy $\textrm{Trace}(E[\psi_i u]; U_i \cap \partial \overline{\Om}^\mathsf{c} ) = \textrm{Trace}(\psi_i u; U_i \cap \partial \Om)$ and $\textrm{Trace}(\nabla E[\psi_i u]; U_i \cap \partial \overline{\Om}^\mathsf{c} ) = \textrm{Trace}(\nabla (\psi_i u); U_i \cap \partial \Om)$.

\

Define the function $E[u]$ via 

$$E[u] := \sum_{i=0}^k E[\psi_i u], $$
where, for the sake of notation, $E[\psi_0 u]$ is just the function $\psi_0 u$ extended by 0 to $\Rn$. As $E[u]$ is the sum of functions in $BH( \Rn ; \Rd)$, it is clearly in $BH(\Rn; \Rd)$, and inside $\Om$ we have

$$E[u] = \sum_{i=0}^k E[\psi_i u] = \sum_{i=0}^k \psi_i u = u.$$
It suffices to verify that $\nabla E[u]$ has the correct trace on $\partial \overline{\Om}^\mathsf{c}$. To see this, note that

\begin{align*}
\textrm{Trace}(\nabla E[u] ; \partial \overline{\Om}^\mathsf{c}) &= \sum_{i=0}^k \ \textrm{Trace}(\nabla E[\psi_i u] ;  \partial \overline{\Om}^\mathsf{c}) \\
&= \sum_{i=0}^k  \ \textrm{Trace}(\nabla(\psi_i u); \partial \Om \cap U_i) \\
&=  \sum_{i=0}^k \ \textrm{Trace}(u \otimes \nabla \psi_i + \psi_i \nabla u ; \partial \Om \cap U_i)\\
&= \sum_{i=0}^k \ \textrm{Trace}(u ; \partial \Om \cap U_i) \otimes \nabla \psi_i +\psi_i \textrm{Trace}(\nabla u ; \partial \Om \cap U_i) \\ 
&= \textrm{Trace}(\nabla u; \partial \Om),
\end{align*}
where in the last line we use the fact that $\sum_{i = 0}^k \nabla \psi_i = \nabla ( \sum_{i=0}^k \psi_i) = \nabla( 1) = 0$.

Since Trace$(\nabla E[u]; \partial \overline{\Om}^\mathsf{c}) = \textrm{Trace}(\nabla u; \partial \Om)$, we conclude that $|D( \nabla E[u])|(\partial \Om) = 0$. \qed

\

We now present the second order version of Theorem \ref{density}.

\bc{appprox} Let $\Om \subset \Rn$ be a bounded Lipschitz domain. For any function $u \in BH(\Om ; \Rd)$ there exist smooth functions $u_n$ such that $u_n \rga u$ in $L^1$, $\nabla u_n \rga \nabla u$ in $L^1$, $\nabla^2 u_n \Ln\wlims D(\nabla u)$, and

$$\int_{\Om} \sqrt{1 + |\nabla^2 u_n|^2} dx  \rga \int_{\Om} \sqrt{1 + |\nabla^2 u|^2} + |D_s (\nabla u)|(\Om).$$
\ec

\proof

Since $\Om$ is Lipschitz, by Theorem \ref{densityL} there is a function $E[u] \in BH(\Rn ; \Rd)$ with $E[u] = u$ in $\Om$ and $|D( \nabla E[u])|(\partial \Om) = 0$.

\

Let $u_n := E[u] * \phi_{1/n}$. Since $\nabla^2 u_n = D(\nabla E[u]) * \phi_{1/n}$, we can apply Theorem \ref{density} to the measure $\mu := D(\nabla E[u])$ using the integrand

$$g(p) := \sqrt{1+|p|^2}, \ g^{\infty}(p) = |p|, $$
noting that $g$ satisifes conditions (\ref{A1}) and (\ref{A3}). 

\qed

\noindent\makebox[\linewidth]{\rule{\paperwidth}{0.4pt}}

\section{Main Result}

In what follows, $\Om \subset \Rn$ is an open, bounded Lipschitz domain. We define the functional $\cF: W^{2,1}(\Om; \Rd) \rga [0, \infty]$ by

$$F[u] := \int_{\Om} f(x, \nabla^2u(x)) dx, \textrm{ for } u \in W^{2,1}(\Om; \Rd), $$
where $f$ is a continuous integrand satisfying the following hypotheses:

\blist
\item[(H1)] Linear growth: $0 \leq f(x,H) \leq C (1 + |H|)$ for all $x \in \Om$, $H \in \mathbb{R}^{d \times N \times N}$ and some $C>0$; 

\item[(H2)] Modulus of continuity: $|f(x, H) - f(y,H)| \leq \omega(|x-y|)(1 + |H|)$ for all $x,y \in \Om, H \in \mathbb{R}^{d \times N \times N}$, where $\omega(s)$ is a nondecreasing function with $w(s) \rga 0$ as $s \rga 0^+$.


\elist
The relaxation of $F$ onto the space $\BHO$ is given by
\begin{align*}
\cF[u] := \inf \bigg\{ \liminf_{n \rga \infty} F[u_n] : \ u_n \rga u \textrm{ in } L^1(\Om; \Rd), \  \sup_n \|u_n\|_{W^{2,1}} < \infty \bigg \}.
\end{align*}

Our main result is the following integral representation theorem.

\bt{mainresult}
\label{mainresult}
If $f$ satisfies (H1) and (H2), then for every $u \in \BHO$ we have

$$\cF[u] =  \int_\Om \cQ_2f(x, \nabla^2u) dx + \int_\Om (\cQ_2f)^{\infty}\bigg(x, \frac{d D_s(\nabla u)}{d |D_s(\nabla u)|} \bigg) d |D_s(\nabla u) |.$$

\et
We will prove this in two steps. Setting
$$\cG[u] := \int_\Om \cQ_2f(x, \nabla^2u) dx + \int_\Om (\cQ_2f)^{\infty}\bigg(x, \frac{d D_s(\nabla u)}{d |D_s(\nabla u)|} \bigg) d |D_s(\nabla u) |, $$
we will show that $\cF \leq \cG$ and $\cG \leq \cF$.

\bt{upperbound}

For all $u \in \BHO$, we have $\cF[u] \leq \cG[u]$.

\et

\proof

We first prove this upper bound for $u \in W^{2,1}(\Om,\Rd)$. By the definition of $\cF$, it suffices to find a sequence of functions $\{ u_n \} \subset W^{2,1}(\Om,\Rd)$ such that $u_n \rga u$ in $L^1(\Om; \mathbb{R}^d)$, $\| u_n \|_{W^{2,1}}$ bounded, and

$$\liminf_{n \rga \infty} \int_\Om f(x, \nabla^2 u_n) dx \leq \int_\Om \cQ_2 f(x, \nabla^2u) dx.$$
The existence of such a sequence is guaranteed by the integral representation of the weakly lower semi-continuous envelope in $W^{2,1}(\Om,\Rd)$ from \cite{BFL00}, Theorem 1.3. In addition, necessary and sufficient conditions on lower-semicontinuity of second order vector valued functionals can be found in \cite{M65}, Theorem 4. Thus, for any $u \in W^{2,1}(\Om, \Rd)$ we have
\begin{equation}
\label{MainUpperBound}
\cF[u] = \cG[u].
\end{equation}

Next we show that for any $u \in \BHO$ we have $\cF[u] \leq \cG[u]$. By Corollary \ref{appprox}, we can find $u_n \in W^{2,1}(\Om,\Rd)$ so that $u_n \rga u$ in $W^{1,1}$, $\nabla^2 u_n \Ln \wlims D(\nabla u),$ and the convergence is area-strict, i.e.,

$$ \lim_{n \rga \infty} \int_{\Om} \sqrt{1 + |\nabla^2u_n|^2} \, dx = \int_{\Om} \sqrt{1 + |\nabla^2u|^2} \, dx + |D_s (\nabla u)|(\Om) .$$
In particular we note that since $u_n \rga u$ in $W^{1,1}$ and $\nabla^2 u_n \Ln \wlims D(\nabla u)$ we must have $\| u_n \|_{W^{2,1}} $ bounded.

Since $\cQ_2 f$ is continuous and 2-quasiconvex with linear growth, Theorem \ref{Area-Strict Continuity} applies. Thus, $\cG$ is continuous with respect to area-strict convergence, and so
$$\cF[u] \leq \liminf_{n \rga \infty} \cF[u_n] = \lim_{n \rga \infty} \cG[u_n] = \cG[u]$$
where we use (\ref{MainUpperBound}) on each of the $u_n$.

\qed

First, we will prove the lower bound for coercive integrands. We claim that in fact, if we have the lower bound for coercive integrands, we have it in general.

\bl{coercive}

If $\cF \geq \cG$ for every integrand satisfying (H1), (H2), and (H3), then $\cF \geq \cG$ for every integrand satisfying (H1) and (H2).

\el

\proof

Let $f$ be a continuous integrand satisfying (H1) and (H2), and consider the coercive integrand $f_{\eps} := f+ \eps|\cdot|$. We observe that

$$\cQ_2(f_{\eps}) \geq \cQ_2f + \eps |\cdot| $$
since $|\cdot|$ is convex. Furthermore, by basic properties of limits,

$$(\cQ_2(f_{\eps}))^{\infty} \geq (\cQ_2f + \eps|\cdot|)^{\infty} = (\cQ_2f)^{\infty} + \eps|\cdot|. $$
Now, for any sequence $\{u_n\} \subset W^{2,1}(\Om,\Rd)$ with $u_n \rga u$ in $L^1$, $\sup_n \|u_n\|_{W^{2,1}} < \infty$, we have
\begin{align*}
\varliminf_{n \rga \infty} &\int_{\Om} f(x, \nabla^2u_n(x)) dx \geq \varliminf_{n \rga \infty} \int_{\Om} \bigg[ f(x, \nabla^2u_n(x)) + \eps |\nabla^2 u_n(x)| \bigg] dx - \varlimsup_{n \rga \infty} \int_{\Om} \eps |\nabla^2u_n(x)| dx \\
&\geq \int_{\Om} \cQ_2f(x, \nabla^2u(x)) dx + \int_{\Om} (\cQ_2f)^{\infty}\bigg(x, \frac{d D_s(\nabla u(x))}{d|D_s(\nabla u(x))|} \bigg) d|D_s(\nabla u(x))| - \eps C
\end{align*}
where we have used the lower bound for $f_{\eps}$ and the fact that $\{\int_{\Om} |\nabla^2u_n| \}$ is bounded. Letting $\eps \rga 0$, we have
$$\varliminf_{n \rga \infty} F[u_n] \geq \cG[u] $$
and, taking the infimum over all such sequences, we conclude that
$$\cF[u] \geq \cG[u].$$
\qed

We will now prove our theorem in the case where $f$ is coercive.

\bt{lowerbound}

Assume that $f$ satisfies (H1), (H2) and (H3). For all $u \in \BHO$, we have $\cG[u] \leq \cF[u]$.

\et

\proof

Let $u \in \BHO$ be given, and let $\{ u_n \} \subset W^{2,1}(\Om,\Rd)$ be an arbitrary sequence with $u_n \rga u$ in $L^1(\Om; \mathbb{R}^d)$ and $\sup_n \| u_n\|_{W^{2,1}} < \infty$. By compactness, we can extract a subsequence (not relabeled) such that $u_n \rga u$ in $W^{1,1}$, $\nabla^2 u_n \ \Ln  \mres \Om \wlims D(\nabla u)$ and the value of

$$\liminf_{n \rga \infty} f(x, \nabla^2 u_n)$$
is unchanged. We proceed according to the blow-up method. Define nonnegative Radon measures $\mu_n$ via

$$
\mu_n(E) := \int_E f(x, \nabla^2 u_n) dx \textrm{ for every Borel set } E \subset \Om.
$$
Without loss of generality we may assume that $\{ \mu_n(\Om) \}$ is bounded, and so, passing to a subsequence (not relabled), we can find a Radon measure $\mu$ such that $\mu_n \wlims \mu$.

We consider the Radon-Nikodym decomposition of $\mu$ with respect to $|D(\nabla u)|$,

$$\mu = \frac{d \mu}{d \Ln} \Ln \mres \Om + \frac{d \mu}{d |D_s(\nabla u)|} |D_s(\nabla u)| + \mu_s,$$
where $\mu_s$ is a nonnegative Radon measure such that $\mu_s \perp |D(\nabla u)|$.

We claim that

\be{Leb}
\frac{d \mu}{d \Ln} (x_0) \geq \cQ_2f(x_0, \nabla^2u(x_0)) \textrm{ for } \Ln a.e. \, x_0 \in \Om,
\ee
and
\be{Sing}
\frac{d \mu}{d |D_s(\nabla u)|} (x_0) \geq (\cQ_2f)^{\infty} \bigg( x_0, \frac{d D_s(\nabla u)}{d |D_s(\nabla u)|}(x_0) \bigg) \textrm{ for } |D_s(\nabla u)|\, a.e. \, \, x_0 \in \Om .
\ee

If $\rf{Leb}$ and $\rf{Sing}$ hold, then we have
\begin{align*}
\liminf_{n \rga \infty} \int_\Om f(x, \nabla^2 u_n) dx &= \liminf_{n \rga \infty} \mu_n(\Om) \geq \mu(\Om) \\
&= \int_\Om \frac{d \mu}{d \Ln} dx + \int_\Om \frac{d \mu}{d |D_s(\nabla u)|} d|D_s(\nabla u)| + \mu_s(\Om) \geq \cG[u].
\end{align*}

The arbitrariness of the sequence $\{u_n\}$ would yield $\cF[u] \geq \cG[u]$. The remainder of this proof is dedicated to proving (\ref{Leb}) and (\ref{Sing}).

\underline{Step 1}: We claim that for $\Ln$ a.e. $x_0 \in \Om$, we have 

$$
\frac{d \mu}{d \Ln} (x_0) \geq \cQ_2f(x_0, \nabla^2u(x_0)).
$$

Note that the measures $\{ |\nabla^2 u_n| \Ln \mres \Om \}$ are bounded in total variation, so, along a subsequence, (not relabled), we have $|\nabla^2 u_n| \wlims \nu$ for some measure $\nu$. By the Lebesgue Differentiation theorem, for $\Ln$ a.e. $x_0 \in \Om$ we have

\be{PwLeb}
\frac{d \mu}{d \Ln} (x_0) = \lim_{\eps \rga 0^+} \frac{ \mu(\overline{\Qeps})}{\eps^N},
\ee
\be{PwLeb2}
\nabla^2u(x_0) = \lim_{\eps \rga 0^+} \frac{D(\nabla u)(\overline{\Qeps})}{\eps^N},
\ee
\be{LebSing}
0 = \lim_{\eps \rga 0^+} \frac{|D_s(\nabla u)|(\overline{\Qeps})}{\eps^N},
\ee
\be{nuLeb}
\frac{d \nu}{d \Ln} (x_0) = \lim_{\eps \rga 0^+} \frac{ \nu(\overline{\Qeps})}{\eps^N} < \infty.
\ee
Select $\eps_k \rga 0$ such that $\mu(\partial Q(x_0, \eps_k)) = \nu(\partial Q(x_0, \eps_k)) = |D(\nabla u)|( \partial Q(x_0, \eps_k)) = 0$, and write 

\begin{align*}
\frac{d \mu}{d \Ln} (x_0) &= \lim_{k \rga \infty} \frac{ \mu(Q(x_0, \eps_k))}{\eps_k^N} = \lim_{k \rga \infty} \lim_{n \rga \infty} \frac{ \mu_n(Q(x_0, \eps_k))}{\eps_k^N} \\
&= \lim_{k \rga \infty} \lim_{n \rga \infty} \frac{1}{\eps_k^N} \int_{Q(x_0, \eps_k)} f(x, \nabla^2 u_n) dx.
\end{align*}
With $x = x_0 + \eps_k y$, we obtain 

\begin{align}
\notag
\frac{d \mu}{d \Ln} (x_0) &= \lim_{k \rga \infty} \lim_{n \rga \infty} \int_Q f(x_0 + \eps_k y, \nabla^2 u_n(x_0 + \eps_k y)) dy\\
&\geq \varlimsup_{k \rga \infty} \varlimsup_{n \rga \infty} \int_Q \cQ_2f(x_0 + \eps_k y, \nabla^2 u_n(x_0 + \eps_k y)) dy.
\label{Q1}
\end{align}
Define functions $v_{n,k} \in W^{2,1}(Q ; \Rd)$ by

$$ v_{n,k}(y) := \frac{u_n(x_0 + \eps_k y) - P_{\eps,k} y - a_{\eps,k} }{\eps_k^2} - \frac{1}{2} \nabla^2 u(x_0)(y,y),$$
where $a_{\eps, k} := \intaverage_Q v_{n,k} $ and $P_{\eps, k} := \intaverage_Q \nabla v_{n,k} $, selected so that each $v_{n,k}$ and its gradient have average zero. By (\ref{Q1}) we get
 
\be{veenkay}
\frac{d \mu}{d \Ln} (x_0) \geq \varlimsup_{k \rga \infty} \varlimsup_{n \rga \infty} \int_Q \cQ_2f(x_0 + \eps_k y, \nabla^2u(x_0) + \nabla^2 v_{n,k}(y)) dy.
\ee

For fixed $k$, the measures $\{ \nabla^2 v_{n,k} \, \Ln \mres Q \}$ converge weakly-$*$ to the measure $\lambda_k$ given by
$$\lambda_k(E) := \frac{D(\nabla u)(x_0 + \eps_k E)}{\eps_k^N} - \nabla^2u(x_0) \Ln(E), \textrm{ for every Borel set } E \subset Q,$$
and by $\{ \lambda_k \}$ converge weakly-$*$ to $0$. to see this, fix any $\psi \in C_c(Q)$. We have

\begin{align*}
\bigg|\int_{Q} \psi(y) d \lambda_k(y)\bigg| &= \frac{1}{\eps_k^N} \bigg|\int_{Q(x_0, \eps_k)} \psi\bigg(\frac{x-x_0}{\eps_k}\bigg) dD(\nabla u)(x) - \int_{Q(x_0, \eps_k)} \psi\bigg(\frac{x-x_0}{r_n}\bigg) \nabla^2u(x_0) dx \bigg|\\
&\leq \frac{1}{\eps_k^N} \bigg| \int_{Q(x_0, \eps_k)} \psi\bigg(\frac{x-x_0}{\eps_k}\bigg) (\nabla^2u(x)-\nabla^2u(x_0)) dx  \bigg|  \\
& \ \ \ + \bigg| \int_{Q(x_0, \eps_k)} \psi\bigg(\frac{x-x_0}{\eps_k}\bigg) dD_s(\nabla u)(x) \bigg|\\
&\leq \|\psi\|_{\infty} \bigg( \frac{1}{\eps_k^N} \int_{Q(x_0, \eps_k)} |\nabla^2u(x) - \nabla^2u(x_0)| + \frac{|D_s(\nabla u)|(Q(x_0, \eps_k))}{\eps_k^N} \bigg)
\end{align*}
which goes to $0$ as $k \rga \infty$ by $(\ref{PwLeb2})$ and $(\ref{LebSing})$. 

We also note that for any $n,k$ we have

$$|\nabla^2 v_{n,k}(y)| \leq |\nabla^2 u_n(x_0 + \eps_k y)| + C $$
for some $C>0$. For fixed $k$ we have

$$|\nabla^2 u_n(x_0 + \eps_k \cdot)| \Ln \mres Q \wlims \frac{T^{\#}_{x_0, \eps_k} \nu}{\eps_k^N}  $$
where $T^{\#}_{x_0, \eps_k} \nu$ denotes the push-forward of $\nu$ under the mapping which takes $x \mapsto \frac{x-x_0}{\eps_k}.$ Since by (\ref{nuLeb}) we have
\begin{align*}
\limsup_{k  \rga \infty } \limsup_{n \rga \infty} \int_{Q} |\nabla^2u_n(x_0 + \eps_k y)| dy &= \limsup_{k \rga \infty} \limsup_{n \rga \infty} \frac{1}{\eps_k^N} \int_{Q(x_0, \eps_k)} |\nabla^2 u_n(x)| dx \\
&\leq \limsup_{k \rga \infty} \frac{\nu(Q(x_0, \eps_k))}{\eps_k^N} < \infty
\end{align*}
we conclude that

$$\limsup_{k  \rga \infty } \limsup_{n \rga \infty} \int_{Q} |\nabla^2 v_{n,k}(y)| dy   < \infty.$$
Thus by Lemma \ref{diags} we can find a diagonalized sequence $v_k := v_{n_k,k}$ such that $\{ \nabla^2 v_k \Ln \mres Q \}$ converges weakly-$*$ to the constant measure $0$ and $\{|\nabla^2 v_k| \Ln \mres Q \}$ converges weakly-$*$ to some nonnegative Radon measure $\pi$. Using the modulus of continuity of $\cQ_2f$, see (\ref{modcontoff}), and (\ref{veenkay}) we have

\begin{align*}
\frac{d \mu}{d \Ln} (x_0) &\geq \lim_{k \rga \infty} \int_Q \cQ_2f(x_0 + \eps_k y, \nabla^2u(x_0) + \nabla^2 v_k(y)) dy \\
&\geq \varlimsup_{k \rga \infty} \int_Q \cQ_2f(x_0, \nabla^2u(x_0) + \nabla^2 v_k(y)) dy  - \int_Q \tilde{C} \omega(\eps_k) (1 + |\nabla^2 v_k(y)| ) dy \\
&\geq \varlimsup_{k \rga \infty}  \int_Q \cQ_2f(x_0, \nabla^2u(x_0) + \nabla^2 v_k(y)) dy
\end{align*}
because $\int_{Q} |\nabla^2 v_k(y)| dy$ is bounded and $\omega(\eps_k) \rga 0$.

In order to apply 2-quasiconvexity, we have to use a $W^{2,1}_0(Q;\mathbb{R}^N)$ perturbation of $\nabla^2u(x_0)$. For $\delta < 1$, let $\phi_{\delta} \in C^{\infty}_c(Q; [0,1])$ be such that $\phi = 1$ on $Q_{\delta} = Q(0, 1- \delta)$, supp$(\phi) \subset Q_{\delta/2} =Q(0,1-\frac{\delta}{2}),$  $\| \nabla \phi \|_{\infty} \leq \frac{C}{\delta}$, $\| \nabla^2 \phi \|_{\infty} \leq \frac{C}{\delta^2}$ for some $C>0$, and let $z_{k, \delta} := \phi_{\delta} \, v_k$. In view of the definition of 2-quasiconvexity (\ref{defQC}), for every $k$ and $\delta$ we have

\begin{equation}
\label{zkdelt}
\int_Q \cQ_2f(x_0, \nabla^2u(x_0) + \nabla^2 z_{k, \delta}(y)) dy \geq \cQ_2f(x_0, \nabla^2u(x_0)).
\end{equation}
\noindent
On the other hand, setting $S_{\delta} := Q_{\delta/2} \setminus \overline{Q_\delta}$, we obtain

\begin{align}
\notag
\int_Q &\cQ_2f(x_0, \nabla^2u(x_0) + \nabla^2 z_{k, \delta}(y)) dy = \int_{Q_{\delta}} \cQ_2f(x_0, \nabla^2u(x_0) + \nabla^2 v_k(y)) dy \\
&+ \int_{S_{\delta}} \cQ_2f(x_0, \nabla^2u(x_0) + \nabla^2 z_{k, \delta}(y))dy + \cQ_2f(x_0, \nabla^2u(x_0)) |Q \setminus Q_{\delta/2} |
\label{strip}
\end{align}
and, as $k$ goes to infinity,

\begin{equation}
\label{cutoff}
\varlimsup_{k \rga \infty} \int_{S_{\delta}} \cQ_2f(x_0, \nabla^2u(x_0) + \nabla^2 z_{k, \delta}(y)) dy \leq
C \lim_{k \rga \infty} \bigg( \delta + \int_{S_{\delta}} \bigg( \frac{1}{\delta^2} |v_k| + \frac{1}{\delta} |\nabla v_k| + |\nabla^2 v_k| \bigg)dy \bigg),
\end{equation}
where we have used the growth condition (H1) and the fact that $\cQ_2f \leq f$.

As we have $\nabla^2 v_k \wlims 0$ and the average of $v_k$ and $\nabla v_k$ are 0, $\{|v_k|\}$ and $\{ |\nabla v_k| \}$ are vanishing in $L^1(Q)$. To see this, let $\{ v_{k_i} \}$ be an arbitrary subsequence of $\{v_k\}$. We observe that, by the Poincar\'e-Wirtinger inequality for $BV$ functions (see \cite{EG92}, 5.6.1), we must have that $\{v_{k_i} \}$ and $\{\nabla v_{k_i}\}$ are bounded in $L^1$. Since we have a bounded sequence in $BH$, we can extract a further subsequence, not relabled, and a function $v \in BH$ such that

$$\lim_{i \rga \infty} \int_Q |v_{k_i}(x)-v(x)| dx = \lim_{i \rga \infty} \int_Q |\nabla v_{k_i}(x) -\nabla v(x)| dx = 0$$
and

$$D(\nabla v_{k_i}) \wlims D(\nabla v) \textrm{ in } Q.$$
However, because $D(\nabla v_{k_i}) \wlims 0$, we have $D(\nabla v) = 0$ and therefore $\nabla v$ is a constant function. Since $Q$ is connected and $\int_{Q} \nabla v = 0$, we must have $\nabla v = 0$. Similarly, this implies that $v$ is a constant function, and $\int_{Q} v = 0$ implies $v = 0$. Thus, we have

$$\lim_{i \rga \infty} \int_Q |v_{k_i}(x)| dx = \lim_{i \rga \infty} \int_{Q} |\nabla v_{k_i}(x)| dx = 0 .$$

Due to the arbitrariness of the subsequence of $\{v_k\}$, we conclude that it is true for our original sequence. As $v_k$ and $\nabla v_k$ are going to 0 in $L^1$, (\ref{cutoff}) becomes

$$  \varlimsup_{k \rga \infty} \int_{S_{\delta}} \cQ_2f(x_0, \nabla^2u(x_0) + \nabla^2 z_{k, \delta}(y)) dy \leq  C \delta + C \pi(\overline{S_{\delta}}).$$
Thus, we have that for every $\delta <1$, using (\ref{zkdelt}) and (\ref{strip}),
\begin{align*}
\varlimsup_{k \rga \infty} \int_{Q} \cQ_2f(x_0, \nabla^2u(x_0) + \nabla^2 v_k(y)) dy &\geq \varlimsup_{k \rga \infty} \int_{Q_{\delta}} \cQ_2f(x_0, \nabla^2u(x_0) + \nabla^2 v_k(y)) dy \\
&\geq \cQ_2f(x_0, \nabla^2u(x_0)) - C \delta - C|Q \setminus Q_{\delta/2}| - C \pi(\overline{S_{\delta}}).
\end{align*}
Note that for every $\delta>0$, $\overline{S_{\delta}} \subset Q \setminus Q_{\delta}$ and $(Q \setminus Q_{\delta}) \searrow \emptyset$ as $\delta \rga 0^+$. Thus, as we let $\delta$ decrease to 0, we have

$$\varlimsup_{k \rga \infty} \int_{Q} \cQ_2f(x_0, \nabla^2u(x_0) + \nabla^2 v_k(y)) dy \geq \cQ_2f(x_0, \nabla^2u(x_0)). $$ 

\

\underline{Step 2}: We show that $|D_s(\nabla u)|$ a.e. $x_0 \in \Om$, we have 

$$
\frac{d \mu}{d |D_s(\nabla u)|} (x_0) \geq (\cQ_2f)^{\infty}\bigg(x_0, \sing(x_0)\bigg).
$$

Since the measures $\{ |\nabla^2 u_n| \Ln \mres \Om \}$ are bounded in total variation, along a subsequence, (not relabled), we have $|\nabla^2 u_n | \wlims \pi$ for some Radon measure $\pi$. Further we can decompose $\pi$ into Radon measures $\pi^A$ and $\pi^B$ such that $\pi^A << |D_s(\nabla u)|$ and $\pi^B \perp |D_s(\nabla u)|$.

By standard properties of BV functions, we know that for $|D_s(\nabla u)|$ a.e. $x_0 \in \Om$ we have

\[ \frac{d \mu}{d |D_s(\nabla u)|} (x_0) = \lim_{\eps \rga 0^+} \frac{ \mu(\overline{\Qeps})}{|D_s(\nabla u)|(\overline{\Qeps})}, \]
\be{PwSing2}
\sing(x_0) = \lim_{\eps \rga 0^+} \frac{D(\nabla u)(\overline{\Qeps})}{|D_s(\nabla u)|(\overline{\Qeps})},
\ee
\be{piSing}
\frac{d \pi^A}{d |D_s(\nabla u)|}(x_0) = \lim_{\eps \rga 0^+} \frac{\pi(\overline{\Qeps})}{|D_s(\nabla u)|(\overline{\Qeps})} < \infty,
\ee
\be{piAC}
\lim_{\eps \rga 0^+} \frac{\pi^B(\overline{\Qeps})}{|D_s(\nabla u)|(\overline{\Qeps})} = 0,
\ee
\be{piACleb}
\lim_{\eps \rga 0^+} \frac{1}{|D_s(\nabla u)|(\overline{\Qeps})} \int_{\Qeps} \bigg| \frac{d \pi^A}{d |D_s(\nabla u)|}(x) - \frac{d \pi^A}{d |D_s(\nabla u)|}(x_0) \bigg| d|D_s(\nabla u)|(x) = 0,
\ee
and
\be{PwSing3}
\lim_{\eps \rga 0^+} \frac{|D_s(\nabla u)|(\overline{\Qeps})}{\eps^N} = \infty.
\ee

Fix $\sigma \in (0,1)$. By Lemma \ref{smallercubes2} we can select $\eps_k \rga 0$ such that $|D_s(\nabla u)|(\partial Q(x_0, \eps_k)) = \mu(\partial Q(x_0, \eps_k)) = \pi(\partial Q(x_0, \eps_k)) = 0$, and

\be{goodeps}
\lim_{k \rga \infty} \frac{|D_s(\nabla u)|(Q(x_0, \sigma \eps_k))}{|D_s(\nabla u)|(Q(x_0, \eps_k))} \geq \sigma^N .
\ee{}
We have

\begin{align*}
\frac{d \mu}{d |D_s(\nabla u)|} (x_0) &= \lim_{k \rga \infty} \frac{ \mu(Q(x_0, \eps_k))}{|D_s(\nabla u)|(Q(x_0, \eps_k))} = \lim_{k \rga \infty} \lim_{n \rga \infty} \frac{ \mu_n(Q(x_0, \eps_k))}{|D_s(\nabla u)|(Q(x_0, \eps_k))} \\
&= \lim_{k \rga \infty} \lim_{n \rga \infty} \frac{1}{|D_s(\nabla u)|(Q(x_0, \eps_k))} \int_{Q(x_0, \eps_k)} f(x, \nabla^2 u_n) dx.
\end{align*}
With the change of variables $x = x_0 + \eps_k y$, we obtain 

\begin{align}
\notag
\frac{d \mu}{d |D_s(\nabla u)|} (x_0) &= \lim_{k \rga \infty} \lim_{n \rga \infty} \frac{\eps_k^N}{|D_s(\nabla u)|(Q(x_0, \eps_k))}  \int_Q f(x_0 + \eps_k y, \nabla^2 u_n(x_0 + \eps_k y)) dy \\
&\geq \varlimsup_{k \rga \infty} \varlimsup_{n \rga \infty} \frac{\eps_k^N}{|D_s(\nabla u)|(Q(x_0, \eps_k))}  \int_Q \cQ_2f(x_0 + \eps_k y, \nabla^2 u_n(x_0 + \eps_k y)) dy.
\label{singQ1}
\end{align}
Note that by (\ref{PwSing3})  

$$t_k := \eps_k^{-N} |D_s(\nabla u)|(Q(x_0, \eps_k)) \rga \infty$$
\noindent
as $k \rga \infty$. We define functions $V_{n,k} \in L^1(Q, \mathbb{R}^{d \times N \times N})$ by

$$V_{n,k}(y) := \frac{1}{t_k} \nabla^2 u_n(x_0 + \eps_k y)  $$
\noindent
and consider the associated matrix-valued measures $\Sigma_{n,k}$

$$\Sigma_{n,k}(E) := \int_{E} V_{n,k}(y) dy = \frac{1}{|D_s(\nabla u)|(Q(x_0, \eps_k))} \int_{x_0 + \eps_k E} \nabla^2 u_n(x) dx $$
\noindent
for every Borel set $E \subset Q$. Note that the total variation of $\Sigma_{n,k}$ is given by

$$|\Sigma_{n,k}|(E)  = \frac{1}{|D_s(\nabla u)|(Q(x_0, \eps_k))} \int_{x_0 + \eps_k E} |\nabla^2 u_n(x)| dx = \frac{(|\nabla^2 u_n| \Ln \mres \Om)(x_0 + \eps_k E) }{|D_s(\nabla u)|(Q(x_0, \eps_k))}.$$
Observe that we can now write (\ref{singQ1}) as

\be{Q1v}
\frac{d \mu}{d |D_s(\nabla u)|} (x_0) \geq \varlimsup_{k \rga \infty} \varlimsup_{n \rga \infty} \frac{1}{t_k} \int_Q \cQ_2f(x_0 + \eps_k y, t_k V_{n,k}(y) ) dy.
\ee
For fixed $k$, we have $\Sigma_{n,k} \wlims \Sigma_k$ and $|\Sigma_{n,k}| \wlims \pi_k$ as $n \rga \infty$, where

\[ \Sigma_k(E) := \frac{D(\nabla u)(x_0 + \eps_k E)}{|D_s(\nabla u)|(Q(x_0, \eps_k))} \textrm{ and } \pi_k(E) := \frac{\pi(x_0 + \eps_k E)}{|D_s(\nabla u)|(Q(x_0, \eps_k))} \]
for every Borel set  $E \subset Q$. Letting $k \rga \infty$, by (\ref{PwSing2}) we have 
\be{SigmaRho}
\Sigma_k \wlims \frac{dD_s(\nabla u)}{d|D_s(\nabla u)|}(x_0) \rho
\ee
where $\rho$ denotes the weak-$*$ limit of $|\Sigma_k|$. This follows from an identical argument to the claim in (\ref{normal}). Indeed, in what follows, we will denote
$$H(x) := \frac{dD_s(\nabla u)}{d|D_s(\nabla u)|}(x) $$
for the sake of notation. Fix $\psi \in C_c(Q)$. We have

\begin{align*}
 \bigg| \int_{Q} \psi(y) d \Sigma_k(y) - \int_{Q} \psi(y) H(x_0) &d |\Sigma_k|(y) \bigg| \\
 &=\frac{1}{|D_s(\nabla u)|(Q(x_0, \eps_k))} \bigg| \int_{Q(x_0, \eps_k)} \psi\bigg(\frac{x-x_0}{\eps_k}\bigg) dD(\nabla u)(x)  \\
& \ \ \ -\int_{Q(x_0, \eps_k)} \psi\bigg(\frac{x-x_0}{\eps_k}\bigg) H(x_0) d |D(\nabla u)|(x) \bigg| \\
&\leq \frac{1}{|D_s(\nabla u)|(Q(x_0, \eps_k))} \bigg| \int_{Q(x_0, \eps_k)} \psi\bigg(\frac{x-x_0}{\eps_k}\bigg) H(x) d |D_s(\nabla u)|(x) \\
& \ \ \ - \int_{Q(x_0, \eps_k)} \psi\bigg(\frac{x-x_0}{\eps_k}\bigg) H(x_0) d |D_s(\nabla u)|(x) \bigg|\\
& \ \ \ +\frac{1}{|D_s(\nabla u)|(Q(x_0, \eps_k))} \bigg( \int_{Q(x_0, \eps_k)} |\psi|\bigg(\frac{x-x_0}{\eps_k}\bigg) |\nabla^2 u|(x) dx \\
& \ \ \ + \int_{Q(x_0, \eps_k)} |\psi|\bigg(\frac{x-x_0}{\eps_k}\bigg) H(x_0) |\nabla^2 u|(x) dx \bigg)\\
&\leq \frac{\| \psi \|_{\infty}}{|D_s(\nabla u)|(Q(x_0, \eps_k))}   \int_{Q(x_0, \eps_k)} |H(x) - H(x_0)| d|D_s(\nabla u)|(x) \\
& \ \ \ + \frac{2 \| \psi \|_{\infty}}{|D_s(\nabla u)|(Q(x_0, \eps_k))} \int_{Q(x_0, \eps_k)} |\nabla^2u(x)|dx,
\end{align*}
which goes to 0 as $k \rga \infty$ in view of (\ref{PwSing2}) and (\ref{PwSing3}). Since

$$\int_{Q} \psi H(x_0)  d|\Sigma_k| = H(x_0) \int_{Q} \psi d |\Sigma_k| \rga H(x_0) \int_{Q} \psi d \rho,$$
we have shown that
$$\int_{Q} \psi d \Sigma_k \rga H(x_0) \int_{Q} \psi d \rho. $$

Note that by (\ref{goodeps}),
\be{fatrho}
\rho(Q) \geq \rho( \sigma \overline{Q}) \geq \varlimsup_k |\Sigma_k|(\sigma \overline{Q}) \geq \varlimsup_k \frac{|D(\nabla u)|(Q(x_0, \sigma \eps_k))}{|D_s(\nabla u)|(Q(x_0, \eps_k))} \geq \sigma^N.
\ee

We also have for any Borel set $E \subset Q$,

\begin{align*}
\notag
\varlimsup_{k \rga \infty} \pi_k(E)  &= \varlimsup_{k \rga \infty}  \frac{\pi(x_0 + \eps_k E)}{|D_s(\nabla u)|(Q(x_0, \eps_k))} =  \varlimsup_{k \rga \infty} \frac{\pi^A(x_0 + \eps_k E)}{|D_s(\nabla u)|(Q(x_0, \eps_k))} \\
\notag
&\leq \varlimsup_{k \rga \infty} \frac{d \pi^A}{d|D_s(\nabla u)|}(x_0) \frac{ |D_s(\nabla u)|(x_0 + \eps_k E)}{|D_s(\nabla u)|(Q(x_0, \eps_k))} \\
\notag
&+ \varlimsup_{k \rga \infty} \frac{1}{|D_s(\nabla u)|(Q(x_0, \eps_k))} \int_{x_0 + \eps_kE} \bigg| \frac{d \pi^A}{d|D_s(\nabla u)|}(x) - \frac{d \pi^A}{d|D_s(\nabla u)|}(x_0) \bigg| d|D_s(\nabla u)|(x) \\
& = \varlimsup_{k \rga \infty} \frac{d \pi^A}{d|D_s(\nabla u)|}(x_0) \frac{ |D_s(\nabla u)|(x_0 + \eps_k E)}{|D_s(\nabla u)|(Q(x_0, \eps_k))}
\end{align*}
by (\ref{piSing}), (\ref{piAC}) and (\ref{piACleb}). Thus, taking $E := Q \setminus \sigma Q$ we see by (\ref{goodeps}) that
\be{sigsok}
\varlimsup_{k \rga \infty} \pi_k(Q \setminus \sigma Q) \leq C (1-\sigma^N)
\ee
for some $C>0$. 
By (\ref{piSing}), we have
\begin{align*}
\varlimsup_{k \rga \infty} \varlimsup_{n \rga \infty} |\Sigma_{n,k}|(Q) &= \varlimsup_{k \rga \infty} \frac{1}{|D_s(\nabla u)|(Q(x_0, \eps_k)} \varlimsup_{n \rga \infty} \int_{Q(x_0, \eps_k)} |\nabla^2u_n(x)| dx  \\
& \leq \varlimsup_{k \rga \infty} \frac{\pi(Q(x_0, \eps_k))}{|D_s(\nabla u)|(Q(x_0, \eps_k))} = \frac{ d \pi^A}{d |D_s(\nabla u)|}(x_0) < \infty
\end{align*}
and by Lemma \ref{diags} we can consider a diagonalized sequence of the $\{ V_{n,k} \}$ and $\{ \Sigma_{n,k} \}$ so that, using (\ref{Q1v}), (\ref{SigmaRho}) and (\ref{fatrho}), we obtain $V_{n_k,k} \Ln \mres Q = \Sigma_{n_k, k} \wlims H_0 \rho$. Note also that
\begin{align*}\varlimsup_{k \rga \infty} \varlimsup_{n \rga \infty} |\Sigma_{n, k}|(Q \setminus \sigma Q) &= \varlimsup_{k \rga \infty} \varlimsup_{n \rga \infty} \frac{1}{|D_s(\nabla u)|(Q(x_0, \eps_k))} \int_{Q(x_0, \eps_k) \setminus Q(x_0, \sigma \eps_k)} |\nabla^2u_n(x)| dx \\
& \leq \varlimsup_{k \rga \infty} \frac{1}{|D_s(\nabla u)|(Q(x_0, \eps_k))} \pi(Q(x_0, \eps_k) \setminus Q(x_0, \sigma \eps_k))  \\
& = \varlimsup_{k \rga \infty} \pi_k(Q \setminus \sigma Q) \leq C (1- \sigma^N)
\end{align*}
by (\ref{sigsok}). Hence, without loss of generality we can choose our diagonalized sequence so that
\be{fatSig}
\varlimsup_{k} |\Sigma_{n_k, k}|(Q \setminus \sigma Q) \leq C (1-\sigma^N)
\ee
 and 

\begin{align*}
\frac{d \mu}{d|D_s(\nabla u)|}(x_0) &\geq \varlimsup_{k \rga \infty} \varlimsup_{k \rga \infty} \frac{1}{t_k} \int_Q \cQ_2f(x_0 + \eps_k y, t_k V_{n,k}(y) ) dy \\
&= \varlimsup_{k \rga \infty} \frac{1}{t_k} \int_Q \cQ_2f(x_0 + \eps_k y, t_k V_{n_k,k}(y) ) dy,
\end{align*}
where
\be{Alberti}
H_0 = H(x_0) = \frac{dD_s(\nabla u)}{d|D_s(\nabla u)|}(x_0) \in \Lambda(N,d,2).
\ee
by the generalized form of the Alberti rank-one theorem (\ref{albertiR1}). Applying the modified modulus of continuity from (\ref{modcontoff}) to (\ref{Q1v}), we have

\begin{align}
\notag
\frac{d \mu}{d |D_s(\nabla u)|} (x_0) &\geq \varlimsup_{k \rga \infty} \int_Q \frac{1}{t_k} \cQ_2f(x_0 + \eps_k y, t_k V_{n_k,k}(y)) dy \\
\notag
&\geq  \varlimsup_{k \rga \infty} \int_Q \frac{1}{t_k} \cQ_2f(x_0 , t_k V_{n_k,k}(y)) dy  - \frac{1}{t_k} \int_{Q} \omega(\eps_k ) ( 1 + |t_k V_{n_k,k}(y)|) dy\\
\notag
&=  \varlimsup_{k \rga \infty} \int_Q \frac{1}{t_k} \cQ_2f(x_0 , t_k V_{n_k,k}(y)) dy  -  \int_{Q} \omega(\eps_k ) \bigg( \frac{1}{t_k} + | V_{n_k,k}(y)| \bigg) dy \\
\label{vKay}
&=  \varlimsup_{k \rga \infty} \int_Q \frac{1}{t_k} \cQ_2f(x_0 , t_k V_{n_k,k}(y)) dy,
\end{align}
\noindent
since $\{\int_{Q} |V_{n_k,k}(y)| dy\}$ is bounded (as $V_{n_k,k} \Ln \mres Q$  is weakly$^*$ converging) and $t_k \rga \infty$. 

We claim that for any $\eta > 0$ we can find $M$ such that if $t > M$ then
\be{poundbound}
\frac{\cQ_2f(x_0, t H)}{t} \geq (\cQ_2f)_{\#}(x_0, H) - \eta
\ee
for every $H$ with $|H| = 1$, where, we recall,
$$ (\cQ_2f)_{\#}(x, H) := \liminf \bigg\{ \frac{\cQ_2f(x',tH')}{t} : x' \rga x, H' \rga H, t \rga \infty \bigg\},$$

To see this, we first prove that the function $(\cQ_2f)_{\#}$ is Lipschitz in its second argument. Indeed, let $x \in \mathbb{R}^N$ and $H,F \in \mathbb{R}^{d \times N \times N}$ be fixed. Let us take any family of sequences $x_n \rga x$, $t_n \rga t$, and $H_n \rga H$. Note that the sequence $F_n := H_n + F - H$ is then converging to $F$, so
$$ (\cQ_2f)_{\#}(x, F)  \leq \varliminf_{n \rga \infty} \frac{\cQ_2f(x_n,t_n F_n)}{t_n} \leq \varliminf_{n \rga \infty} \frac{\cQ_2f(x_n,t_n H_n)}{t_n}  + L |F-H|$$
where $L$ is the Lipschitz constant of $\cQ_2f$ from Lemma \ref{LipCon}, noting that this constant depends only on the growth constant $C$ and the dimensions $N$ and $d$. But, since this holds for any $x_n \rga x$, $t_n \rga t$, $H_n \rga H$, we must have
$$(\cQ_2f)_{\#}(x, F)  \leq (\cQ_2f)_{\#}(x, H)  + L |F-H|  $$
and, by symmetry, we conclude that $(\cQ_2f)_{\#}$ is Lipschitz continuous in its second argument. Now, we claim that (\ref{poundbound}) holds. If not, then there exists $\eta >0$, $\{H_n\}$ with $|H_n| = 1$ and $t_n \rga \infty$ such that 
\be{contradict}
\frac{\cQ_2f(x_0, t_n H_n)}{t_n} < (\cQ_2f)_{\#}(x_0, H_n) - \eta.
\ee
Without loss of generality, since the unit sphere is compact, we can assume that $H_n \rga H$ for some $H$ with $|H| = 1$. Thus, letting $n \rga \infty$ in (\ref{contradict}), we obtain
\begin{align*}
\cQ_2f_{\#}(x_0, H) \leq \varliminf_{n \rga \infty} \frac{\cQ_2f(x_0, t_n H_n)}{t_n} &\leq \lim_{n \rga \infty} \cQ_2f_{\#}(x_0, H_n) - \eta \\
& = \cQ_2f_{\#}(x_0, H) - \eta
\end{align*}
where we have used the fact that $(\cQ_2f)_{\#}$ is Lipschitz continuous in its second argument. This is a contradiction, so we conclude (\ref{poundbound}). In turn, (\ref{poundbound}) implies that for any $\tilde{H}$ and $t$ such that $|\tilde{H}| > \frac{M}{t}$ we have
\be{Qpound}
\frac{\cQ_2f(x_0, t \tilde{H})}{t} \geq (\cQ_2f)_{\#}(x_0, \tilde{H}) - \eta |\tilde{H}|
\ee
by letting $H := \frac{\tilde{H}}{|\tilde{H}|}$ in (\ref{poundbound}). Consider the set

$$E_k := \bigg\{ x \in Q : | V_{n_k,k}(x)| > \frac{M}{t_k} \bigg\}.$$
We have by (\ref{Qpound})

\begin{align}
\notag
\int_Q \frac{1}{t_k} \cQ_2f(x_0 , t_k V_{n_k,k}(y)) dy &\geq  \int_{E_k} \frac{1}{t_k} \cQ_2f(x_0 , t_k V_{n_k,k}(y)) dy \\
&\geq \int_{E_k} \bigg[ (\cQ_2f)_{\#}(x_0, V_{n_k,k}(y)) - \eta |V_{n_k,k}(y)| \bigg] dy \nonumber \\
&\geq \int_{E_k} (\cQ_2f)_{\#}(x_0, V_{n_k,k}(y))dy  - \eta \int_Q |V_{n_k,k}(y)| dy \label{eta}
\end{align}
We can write

\begin{align*}
\int_{Q} (\cQ_2f)_{\#}(x_0, V_{n_k,k}(y))dy &= \int_{E_k} (\cQ_2f)_{\#}(x_0, V_{n_k,k}(y))dy + \int_{Q \setminus E_k} (\cQ_2f)_{\#}(x_0, V_{n_k,k}(y))dy\\
&\leq \int_{E_k} (\cQ_2f)_{\#}(x_0, V_{n_k,k}(y))dy + \int_{Q \setminus E_k} C |V_{n_k,k}(y)| dy \\
&\leq \int_{E_k} (\cQ_2f)_{\#}(x_0, V_{n_k,k}(y))dy + C \frac{M}{t_k},
\end{align*}
and thus

\be{elikecube}
\int_{E_k} (\cQ_2f)_{\#}(x_0, V_{n_k,k}(y))dy \geq \int_{Q} (\cQ_2f)_{\#}(x_0, V_{n_k,k}(y))dy - C \frac{M}{t_k}.
\ee

As discussed in Theorem \ref{Area-Strict Continuity}, the function $(\cQ_2f)_{\#}(x_0, \cdot)$ is positively 1-homogeneous and $\Lambda(N,d,2)$-convex. Since $H_0 \in \Lambda(N,d,2)$ (see (\ref{Alberti})), by \cite{KK16}, Theorem 1.1, we can find an affine function $L(H) = b + \xi \cdot H$ such that $L \leq (\cQ_2f)_{\#}(x_0, \cdot)$ and

\be{contact}
L(H_0) =  (\cQ_2f)_{\#}(x_0, H_0).
\ee
Then, 

\begin{align}
\notag 
\limsup_{k \rga \infty} \int_Q (\cQ_2f)_{\#}(x_0, V_{n_k,k}(y)) dy &\geq  \limsup_{k \rga \infty} \int_Q \Big( b + \xi \cdot V_{n_k,k}(y) \Big) dy \\
\label{linear}
& = b + \limsup_{k \rga \infty} \int_Q \xi \cdot d \Sigma_{n_k,k}(y). 
\end{align}
Let $\psi_{\sigma} \in C_c(Q; [0,1])$ be such that $\psi_{\sigma} = 1$ in $\sigma Q$. We have

\begin{align}
\notag
\limsup_{k \rga \infty} \int_Q \xi \cdot d \Sigma_{n_k,k}(y) &= \lim_{k \rga \infty} \int_Q \psi_{\sigma}(y) \xi \cdot d \Sigma_{n_k,k}(y) +  \limsup_{k \rga \infty} \int_Q (1- \psi_{\sigma}(y)) \xi \cdot d \Sigma_{n_k,k}(y) \\
\notag
&= \int_Q \xi \cdot H_0 \psi_{\sigma}(y) d \rho(y)  +  \limsup_{k \rga \infty} \int_Q (1- \psi_{\sigma}(y)) \xi \cdot d \Sigma_{n_k,k}(y)\\
\label{fatsig}
& \geq \xi \cdot H_0  \ \sigma^N -  C \limsup_{k \rga \infty} |\Sigma_{n_k, k}|(Q \setminus \sigma Q)
\end{align}
for some $C > 0$, where we used (\ref{fatrho}).
In view of (\ref{fatSig}),
\be{thinsig}
\limsup_{k \rga \infty} |\Sigma_{n_k, k}|(Q \setminus \sigma Q) \leq C(1- \sigma^N),
\ee
therefore, by (\ref{linear}),
\[ \limsup_{k \rga \infty} \int_Q \xi \cdot d \Sigma_{n_k,k}(y) \geq  \xi \cdot H_0  \sigma^N  - C(1- \sigma^N), \]
and so by (\ref{linear}), (\ref{fatsig}), and (\ref{thinsig}) 
\begin{align*}
\liminf_{k \rga \infty} \int_Q (\cQ_2f)_{\#}(x_0, v_{n_k,k}(y)) dy &\geq b + \xi \cdot H_0 \sigma^N - C(1-\sigma^N)\\
&\geq  \sigma^N (\cQ_2f)_{\#}(x_0,H_0) - C (1-\sigma^N),
\end{align*}
where we used (\ref{contact}). Putting this together with (\ref{vKay}), (\ref{eta}) and (\ref{elikecube}), we have
\begin{align*}
\frac{d\mu}{d|D_s(\nabla u)|}(x_0) &\geq \varlimsup_{k \rga \infty} \int_Q \frac{1}{t_k} \cQ_2f(x_0 , t_k V_{n_k,k}(y)) dy \\
&\geq \varlimsup_{k \rga \infty} \bigg[ \int_{E_k} (\cQ_2f)_{\#}(x_0, V_{n_k,k}(y))dy  - \eta \int_Q |V_{n_k,k}(y)| dy \bigg] \\
&\geq \limsup_{k \rga \infty} \bigg(  \int_{Q} (\cQ_2f)_{\#}(x_0, V_{n_k,k}(y))dy - C \frac{M}{t_k} \bigg) - \eta C \\
&\geq \sigma^N (\cQ_2f)_{\#}(x_0,H_0) - C (1-\sigma^N) - \eta C.
\end{align*}
Given the arbitrariness of $\eta >0$, we conclude that 

\begin{align*}
\frac{d\mu}{d|D_s(\nabla u)|}(x_0) \geq \sigma^N (\cQ_2f)_{\#}(x_0,H_0) - C (1-\sigma^N)
\end{align*}
for every $\sigma \in (0,1)$, thus, as we send $\sigma \rga 1^{-}$ we conclude that

$$\frac{d\mu}{d|D_s(\nabla u)|}(x_0) \geq (\cQ_2f)_{\#}(x_0,H_0) = (\cQ_2f)^{\infty}(x_0, H_0)$$
since, as noted in (\ref{sameAlb}), $\cQ_2f_{\#} = \cQ_2f^{\infty}$ on $\Lambda(N,d,2)$.

\qed

\section{Acknowledgements}

The author thanks I.\ Fonseca for her invaluable insights and guidance throughout this project. The author wishes to acknowledge the Center for Nonlinear Analysis where this research was accomplished. This work was partially funded by the National Science Foundation under Grant No. DMS-1411646 and NSF PIRE Grant No. OISE-0967140.

\end{document}